\documentclass[opre,nonblindrev]{Revised_informs3}
\SingleSpacedXI

\usepackage{endnotes}
\let\footnote=\endnote
\let\enotesize=\normalsize
\def\notesname{Endnotes}%
\def\makeenmark{\hbox to1.275em{\theenmark.\enskip\hss}}
\def\enoteformat{\rightskip0pt\leftskip0pt\parindent=1.275em
  \leavevmode\llap{\makeenmark}}

\usepackage{natbib}

 \bibpunct[, ]{(}{)}{,}{a}{}{,}%
 \def\bibfont{\small}%
 \def\bibsep{\smallskipamount}%
 \def\bibhang{24pt}%

\usepackage{geometry}
\TheoremsNumberedThrough     % Preferred (Corollary 1, Theorem 1, Theorem 2)
\ECRepeatTheorems

\EquationsNumberedThrough    % Default: (1), (2), ...
\usepackage{amsmath} %% Hossein added this package
\usepackage{wasysym} %% Hossein added this package
\usepackage{algorithm}
\usepackage{algorithmic}

\usepackage{multirow}

\usepackage{pdflscape} %landscape

\usepackage{graphicx} %resize

\usepackage{rotating} %tablesideways

\usepackage{longtable, lscape} %longtable

\newcommand{\specialcell}[2][c]{%
	\begin{tabular}[#1]{@{}c@{}}#2\end{tabular}} %for breaking line in a Table's cell.

\usepackage{mathtools}

\EquationsNumberedThrough    % Default: (1), (2), ...
\usepackage{amsmath} %% Hossein added this package
\usepackage{wasysym} %% Hossein added this package
\usepackage{algorithm}
\usepackage{algorithmic}

\usepackage{multirow}

\usepackage{pdflscape} %landscape

\usepackage{graphicx} %resize

\usepackage{rotating} %tablesideways

\usepackage{longtable, lscape} %longtable

\usepackage{lmodern,bm}  

\usepackage{color,soul} %% for gihlighting: example >>>> \hl{hossein.hashemi@encs.concordia.ca}
\usepackage{bbm} %for \mathbbm{1}

\usepackage[english]{babel} %%% 'french', 'german', 'spanish', 'danish', etc.
\usepackage{amssymb}
\usepackage{amsmath}
\usepackage{mathdots}

\usepackage{tocloft}

\usepackage{etoolbox}

\usepackage{setspace}
\usepackage{enumitem}

\usepackage{afterpage} %To prevent page break before landscape tables (check afterpage command in the text which takes the whole landscapce page as the input)

\usepackage{color,soul} %% for gihlighting: example >>>> \hl{hossein.hashemi@encs.concordia.ca}

\usepackage{booktabs} %for resizing tabular by resiebox

\usepackage{color,soul} %% for gihlighting: example >>>> \hl{hossein.hashemi@encs.concordia.ca}

\usepackage{verbatim} %for deactiviating (commenting) a part of the text. Works with \begin{comment} ...  \end{comment}

\usepackage{afterpage} % \afterpage{ ...} This is for pagebreak for landscape commad

\usepackage{changepage}

\makeatletter
\newcommand\makebig[2]{%
	\@xp\newcommand\@xp*\csname#1\endcsname{\bBigg@{#2}}%
	\@xp\newcommand\@xp*\csname#1l\endcsname{\@xp\mathopen\csname#1\endcsname}%
	\@xp\newcommand\@xp*\csname#1r\endcsname{\@xp\mathclose\csname#1\endcsname}%
}
\makeatother

\makebig{biggg} {3.0}
\makebig{Biggg} {3.5}
\makebig{bigggg}{4.0}
\makebig{Bigggg}{4.5}
\makebig{Biggggg}{5}
\makebig{Bigggggg}{6}
\makebig{Biggggggg}{7}
\makebig{Bigggggggg}{8}

\begin{document}
%%%%%%%%%%%%%%%%
%%\setlength{\footskip}{60pt}
%%\addtolength\footskip{5cm}

\setcounter{tocdepth}{0}

%\setcounter{page}{1}
%\DOI{10.1287/xxxx.0000.0000}%

\RUNAUTHOR{Hashemi Doulabi et al.}

\RUNTITLE{Operating room planning with pooling downstream beds among specialties}

% \TITLE{Bundling Information Goods of Decreasing Value}
\TITLE{Operating room planning with pooling downstream beds among specialties: A stochastic programming approach}
%%\setlength{\TITLE}{-10em}
\ARTICLEAUTHORS{%

		\AUTHOR{Arian Andam$^{1,2}$, Hossein Hashemi Doulabi$^{1,2}$\thanks{Word count (including tables and figures, excluding references and appendices): 8,158 words.}}
\AFF{$^{1}$Department of Mechanical and Industrial Engineering, Concordia University, Montreal, Canada \\$^{2}$Interuniversity Research Center on Enterprise Networks, Logistics and Transportation (CIRRELT), Montreal, Canada\\ \EMAIL{arian.andam@mail.concordia.ca, hossein.hashemi@concordia.ca}} %, \URL{}}
} % end of the block

\ABSTRACT{%
In this paper, we study pooling downstream beds among different 
specialties in a stochastic operating room planning problem. The main sources of uncertainty are stochastic surgical 
durations and patients' lengths of stay. We developed a two-stage 
stochastic programming model where in the first stage we decide on 1) 
the number of non-shared ICU and ward beds to be allocated to each specialty, and 
2) the allocation of surgeries to operating rooms during the planning 
horizon. In the second stage, we decide on 1) the number of shared beds in 
ICU and wards to be allocated to different specialties on each day during the 
planning horizon, 2) the surge capacity required to satisfy downstream 
service to patients, and 3) the overtime incurred in operating rooms. The proposed model aims at minimizing the 
total cost including the patients' waiting cost, postponement cost, 
overtime and fixed cost of operating rooms, and the cost of downstream 
surge capacity. We have implemented the proposed stochastic programming 
model in a sample average approximation framework. To enhance the efficiency of sample average approximation, we have developed a specialized algorithm that quickly solves the second-stage model for any given first-stage
solution for a large number of scenarios. We have carried out extensive computational experiments to evaluate the effectiveness of 
several pooling policies for downstream beds and also the efficiency of 
the proposed sample average approximation algorithm. Moreover, we have 
performed an extensive sensitivity analysis of cost and stochastic 
parameters. Our results demonstrated that a full-sharing policy among 
different specialties in the downstream units enhance the functionality 
of the system by up to 19.53\%. Moreover, the results indicated that the 
solutions obtained by the proposed stochastic model outperform those 
from the corresponding deterministic problem by 17.43\% on average.
}%

\KEYWORDS{Operating room planning, pooling of downstream beds, two-stage 
	stochastic programming, sample average approximation.}
%\HISTORY{Submitte on February 07, 2018; Revised on April 16 and July 23, 2019, Accepted on July 31, 2019}
\maketitle
%%%%%%%%%%%%%%%%%%%%%%%%%%%%%%%%%%%%%%%%%%%%%%%%%%%%%%%%%%%%%%%%%%%%%%
\vspace{-15pt}
\section{Introduction}
Despite all the shocks and abnormalities, which have influenced the healthcare industry in recent years, this sector has successfully managed to keep its essence and structure. In any healthcare system, managers try to keep expenses at a minimum level. In this regard, many health practitioners use mathematical models to identify the bottlenecks of their planning and scheduling systems and fix them (e.g., see \citealp{haghi2023integrated, abbaszadeh2025drone, farghadani2025stochastic}). Based on the existing literature, operating rooms are the most crucial part of hospitals and account for the majority of expenses and revenues. \cite{umali2020efficiency} and \cite{Research2019Global} state that almost two-thirds of the revenues of hospitals come from operating rooms while taking account for 40\% of the expenses.

In this context, the main role of health practitioners is to allocate surgeries to operating rooms over a planning horizon such that medical resources including the available times of operating rooms are used as efficiently as possible. This planning procedure is a very difficult task because the health practitioners must take into account many other restrictive details such as the limited number of beds in ICU and wards. Efficient management of operating rooms considering the limited downstream resources is even a more complicated task because the patients' lengths of stay (LOS) are uncertain and therefore the health practitioners have a hazy view of the available number of beds in ICU and wards within the next few days.

In practice, ICU and ward beds are divided between different specialties to avoid any conflict between the surgical groups. However, a fixed and inflexible allocation of downstream resources to specialties could cause inefficient use of them. For instance, a waste of recourses happens when the patients of one specialty stay shorter than expected in ICU while patients of another specialty require to stay in ICU longer than expected. In this case, the future surgeries of the latter specialty may need to be canceled due to unavailable ICU beds, while extra ICU beds are available for the other surgical group. A question that arises in this context is whether having some shared downstream beds could improve the management of operating rooms and if it does by what margin it enhances the efficiency of the system.

Many papers in the literature have studied operating room planning problems with various assumptions and solution approaches (\citealp{hashemi2014constraint, hashemi2016constraint, hashemi2020vehicle,roshanaei2017collaborative, roshanaei2017propagating, roshanaei2020branch, roshanaei2020reformulation, roshanaei2021solving,breuer2020robust,naderi2021increased,aringhieri2022combining,bargetto2023branch, dolatkhah2026rlcg}). In the context of multi-specialty operating room planning, \cite{castro2015operating} developed a decomposition-based approach for short-term surgical scheduling that coordinates surgical assignments across operating rooms and aligns the schedules of surgeons working in different rooms. In this work, they focused on operating room utilization without explicitly modeling downstream bed capacity. \cite{jung2019scheduling} studied operating room capacity allocation under uncertainty in emergency patient arrivals, and developed optimization and rescheduling procedures to generate weekly and daily schedules. \cite{andam2021operating} developed a mixed-integer programming model for operating room planning that explicitly accounts for multiple downstream units, including wards and ICUs, while minimizing operating room opening, overtime, refusal, and patient waiting costs. \cite{li2024scheduling} studied an elective surgery scheduling problem with integrated resource constraints across preoperative, perioperative, and postoperative stages and proposed a mixed-integer programming model and heuristics. \cite{yin2025multi} studied a multi-specialty operating room scheduling problem with downstream recovery beds at an inter-hospital level, emphasizing on metaheuristic solution approaches. 

On the other hand, many researchers have focused on stochastic optimization approaches such as sample average approximation (SAA) as the solution approach (\citealp{mancilla2012sample, xiao2016stochastic, eun2019scheduling, guo2021logic,coban2023effect,tsang2025stochastic,almoghrabi2025surgery}). \cite{mancilla2012sample} used sample average approximation to model uncertain service durations in a stochastic integer program for the scheduling of single operating room, capturing waiting, idle, and overtime costs through sampled scenarios. \cite{xiao2016stochastic} used SAA within a three-stage recourse framework for OR scheduling with cancellations and resource disruptions, enabling evaluation of overtime and cancellation policies under sampled uncertainty. \cite{eun2019scheduling} applied SAA to approximate a stochastic mixed-integer OR planning model with uncertain surgery durations, balancing overtime and patient health deterioration via scenario-based assignments. \cite{guo2021logic} adopted SAA to approximate uncertain surgery durations in a stochastic distributed operating room scheduling model across hospitals, embedding the scenario-based formulation in decomposition algorithms to obtain robust assignments.  \cite{coban2023effect} applied SAA to single-operating-room outpatient surgery sequencing with uncertain case durations, using sampled distributions to evaluate schedules and compare sequencing heuristics under limited historical data. \cite{tsang2025stochastic} solved a stochastic OR allocation, assignment, sequencing, and scheduling model under duration uncertainty using SAA, where sampled scenarios represent uncertainty across multiple interrelated resource decisions. \cite{almoghrabi2025surgery} benchmarked against SAA-based two-stage stochastic models that sample emergency arrivals and surgery durations to estimate expected cancellation, waiting, and overtime costs in elective–emergency operating room planning. 

Some researchers have focused on the operating room planning problem with the limited number of beds downstream while addressing the stochastic nature of patients’ lengths of stay (\citealp{min2010scheduling, jebali2015stochastic, jebali2017chance,zhang2019two,zhang2020column,zhang2021two, zhan2025optimizing}). \cite{min2010scheduling} proposed a two-stage stochastic programming model for an operating room planning problem with a limited number of ICU beds, where bed capacity is enforced as a hard constraint. \cite{jebali2015stochastic,jebali2017chance} extended this work to limited ward beds, with and without emergency patients. \cite{zhang2019two} developed a similar two-stage stochastic programming model for an operating room planning problem, where patients are denied admission to ICU when all beds are occupied. They proposed an approximate dynamic programming model as the solution approach. \cite{zhang2020column} further proposed a column-generation-based heuristic for a modified version of \cite{min2010scheduling}. \cite{zhang2021two} developed a two-level framework combining Markov decision processes and stochastic programming to account for inter-temporal effects and ICU capacity. In a related direction, \cite{rachuba2022tactical} developed a chance-constrained model to derive admission templates considering operating room and ICU capacities under stochastic demand. \cite{zhan2025optimizing} formulated a two-stage stochastic programming model that jointly considers operating room scheduling and the capacity of post-anesthesia care unit, where recourse costs are evaluated via discrete-event simulation.

There is also another category of papers in the literature that have proposed robust optimization models for operating room planning problems with downstream constraints. \cite{wang2017distributionally} developed a distributionally robust chance-constrained framework to schedule surgeries over a planning horizon considering the restrictions in downstream. They considered uncertain surgical durations in their model but supposed that lengths of stay are deterministic. \cite{neyshabouri2017two} proposed a two-stage robust optimization model to allocate patients to multiple surgery blocks and developed a column-and-row generation algorithm to solve it. In their proposed model, they considered penalties for the violation of the ICU bed constraints in the second stage. \cite{moosavi2018scheduling} formulated a two-stage robust optimization model for an operating room scheduling problem with upstream and downstream constraints. \cite{moosavi2020robust} offered a two-stage heuristic algorithm to solve a scenario-based robust operating room planning problem with uncertain lengths of stay, surgical durations, and emergency demands. \cite{shehadeh2021distributionally} proposed a distributionally robust optimization model to allocate elective surgeries to surgical blocks. They reformulated the problem as a nonlinear-mixed-integer programming model and then linearized it. \cite{augustin2022data} also developed a data-driven-based mixed integer programming model to consider the availability of ICU beds in operating room planning. They calculate the probability of cancellations of surgeries on each day using the historical data. Moreover, \cite{koushan2025sustainable} proposed a scenario-based robust multiobjective operating room scheduling model that incorporates ICU and ward beds while accounting for staff-related sustainability objectives. Furthermore, \cite{shehadeh2025operating} developed a distributionally robust optimization framework for integrated surgery assignment, sequencing, and scheduling with downstream recovery units under LOS ambiguity.

In the literature, limited downstream constraints have also been studied in the context of master problem scheduling in which a cyclic weekly plan must be generated by allocating operating rooms to specialties. \cite{belien2007building} studied a master surgery scheduling problem assuming the number of patients in operating rooms and the patients’ lengths of stay are uncertain. They formulate the problem as a mixed integer programming model minimizing the maximum required number of beds in the objective function. \cite{belien2009decision} developed a decision support system for cyclic master surgery scheduling with the explicit goal of leveling downstream bed occupancy while limiting operating room sharing across specialties and preserving schedule regularity. Their approach combines mixed integer programming with metaheuristic techniques to address multiple, potentially conflicting objectives. \cite{kumar2018sequential} developed a stochastic mixed integer programming model to manage the patient flow in a master surgery scheduling problem with stochastic lengths of stay. In the objective function, they maximized the sum of lengths of stay for scheduled patients in the planning horizon. \cite{makboul2025multiobjective} developed a multi-objective robust master surgical scheduling model that accounts for ICU downstream resources and LOS uncertainty, using simulation-based risk measures to evaluate solution robustness.

Researchers in the literature have addressed the concept of pooling in healthcare systems from different perspectives. \cite{el2023tailored} studied inventory pooling in pharmaceutical supply chains under demand uncertainty and formulated a mixed-integer conic quadratic optimization model to balance pooling costs and safety stock savings. \cite{vanberkel2012efficiency} analyzed the trade-off between centralized and decentralized hospital service structures using queueing models and simulation to identify conditions under which pooling improves efficiency. \cite{batun2011operating} quantified the benefits of pooling operating rooms by developing a two-stage stochastic mixed-integer programming model and exploiting structural properties to solve realistically sized scheduling problems. Moreover, several studies have focused specifically on bed pooling and bed overflow management in hospitals. \cite{lee2018strategies} examined pooling versus focusing of inpatient beds and proposed a two-stage framework combining clustering methods with mixed-integer nonlinear programming for bed allocation decisions. \cite{qiu2025two} studied cross-hospital bed pooling using a two-stage approach integrating weighted K-means clustering with mixed-integer linear programming to minimize social costs and patient waiting times. \cite{izady2021clustered} proposed a clustered overflow bed configuration and developed heuristic optimization models that balance patient rejection and nursing costs. \cite{dai2019inpatient} modeled inpatient bed overflow decisions using a multiclass, multipool queueing framework and applied approximate dynamic programming to derive effective real-time overflow policies under uncertainty. \cite{alsabah2015assigning} addressed inpatient room pooling in Kuwaiti hospitals using a hybrid approach based on queueing theory, data-driven clustering, and integer programming. \cite{hadid2025digital} investigated fully pooled bed allocation across specialties using a digital-twin–enabled simulation-based optimization framework and demonstrated significant improvements in bed utilization and waiting times.

Although there is a rich literature on stochastic operating room planning and healthcare capacity pooling, to the best of our knowledge, there is no paper studying the effect of pooling downstream beds in the presence of uncertain lengths of stay. An earlier version of this work was included in the Master's thesis of the first
author (\citealp{Andam2022OR}). The main contributions of this research are as follows:

\begin{itemize}
	\item For the first time, we study the pooling of downstream beds among different specialties in an operating room planning problem where surgical durations and patients' lengths of stay are stochastic. Our goal is to determine whether the pooling of downstream beds results in enhancing the efficiency of the operating room planning and if so by what margin. 
	\item We formulate the problem as a two-stage stochastic programming model where surgical durations and patients' lengths of stay in downstream resources are stochastic. We embed the proposed stochastic programming model in a sample average approximation (SAA) framework. We also devise a specialized algorithm that enables us to solve the second-stage model for a given first-stage solution with a large number of scenarios. This significantly reduces the standard deviation of the statistical upper bound. As mentioned earlier, there is a rich literature on the use of the SAA algorithm for operating room planning and scheduling problems; however, compared with this literature, this work is the first to apply SAA to evaluate the effectiveness of bed pooling policies for downstream units in an operating room planning context.
	\item We provide extensive computational results to evaluate the improvement in the objective value as a result of sharing beds among specialties for three different pooling policies. We compare our results with the case that a simple deterministic integer programming model is used to solve the problem. We also evaluate the efficiency of the proposed sample average approximation algorithm. Moreover, we perform a broad sensitivity analysis to evaluate the behavior of the proposed stochastic programming approach with different parameter settings. 
\end{itemize}

The remainder of this paper is structured as follows: In Section 2, we have described the problem definition. In Section 3, we have proposed a two-stage stochastic model for the operating room planning problem with the possibility of sharing downstream beds among different specialties. In Section 4, we have proposed a sample average approximation approach to obtain statistical bounds from the developed stochastic programming model. In Section 5, we have carried out extensive computational experiments. Finally, we have concluded this research in Section 6 and have provided some suggestions for future research.

\section{Problem definition}

We consider an operating room planning problem where several patients must be allocated to operating rooms over a planning horizon. In each operating room, there is a limited regular time available for performing surgeries. The durations of surgeries are stochastic and therefore there is a possibility to have overtime at the end of the available time of operating rooms. There is a fixed cost associated with opening each operating room for each day. Patients belong to different surgical specialties that do not share operating rooms. This means that patients allocated to the same operating room must be associated with the same surgical specialty. Each patient has a time window that indicates the earliest and the latest days on which the surgery can be performed. If the latest day is within the planning horizon, we refer to the patient as a mandatory patient that must be operated in the planning horizon before the deadline. Otherwise, the patient is an optional patient that can be postponed to the next planning horizon for the operation. 

After surgeries, the patients are transferred to ICU and wards consecutively and usually stay in each area for a few days. Each surgical specialty has its own number of beds in ICU and wards. For example, cardiovascular, neurology, and orthopedic may have 12, 10, and 9 beds in ICU, respectively. As in the common practice, we suppose that these beds are reserved for the corresponding specialty and will not be occupied by the patients of other specialties in any case. Besides, to study the effect of pooling downstream beds, we suppose that there is a limited number of ICU and ward beds that are shared between specialties. Therefore, in the case that a new cardiovascular patient is operated and needs an ICU bed, but all beds of this specialty are already occupied, the patient may be allocated to one of the available shared ICU beds. If there is not any available ICU bed in the cardiology part of ICU or the shared ICU area, the patient will use a surge capacity at an extra cost. The surge capacity refers either to hiring a part-time ICU nurse, or paying an existing ICU nurse overtime to take care of the patient, or transferring the patient to another hospital with available ICU beds. We consider that there is a similar limited bed capacity for different specialties in wards area with the possibility of having shared beds as well. In this problem, management of downstream beds is specifically more difficult because the patients' lengths of stay are stochastic and therefore the number of patients occupying beds does not follow a clear trend as patients may stay longer or shorter than expected.

In the next section, we propose a two-stage stochastic programming model to formulate this problem. In the first stage, the decision maker determines 1) the number of ICU and ward beds to be allocated to each specialty, and 2) the allocation of surgeries to operating rooms during the planning horizon. Here, we suppose that there is a limited number of beds that can be considered shared between different specialties. This is necessary in practice to avoid conflict between surgical groups. In the second stage, the decision maker determines 1) the number of shared beds in ICU and wards to be allocated to different specialties on each day during the planning horizon, 2) the surge capacity required to satisfy downstream service to patients, and 3) the overtime incurred in operating rooms during the planning horizon. These variables are computed based on first-stage decisions on the allocation of surgeries to operating rooms and also the allocation of non-shared downstream beds to specialties. The objective function is to minimize the total cost including the patients' waiting cost, postponement cost, overtime and fixed cost of operating rooms, and the cost of downstream surge capacity.

\section{Two-stage stochastic programming model}

In this section, we developed a two-stage stochastic programming model. In the first stage, we decide on the allocation of the downstream beds to specialties, which operating rooms to open on each day, and also allocate surgeries to different operating rooms during the planning horizon. In the second stage, we determine how to allocate the shared beds in ICU and wards to different specialties on each day in the planning horizon. We also compute the incurred overtime in operating rooms and the surge capacity in ICU and wards. The first-stage objective function minimizes the patients' waiting cost, postponement cost, and fixed cost of operating rooms, while the second-stage objective function considers the overtime cost of operating rooms and the cost of downstream surge capacity.

\subsection{First-stage model}

The sets, parameters, and variables in our model are as follows. 

\vspace{-5pt}

\noindent
Sets:\vspace{6pt}

	\hspace{-5pt} \begin{tabular}{p{12pt}p{5pt}p{15.45cm}} 
	$\mathcal{I}$ & : & The set of surgeries (or patients). We have $\mathcal{I}=\mathcal{I}_{1}\cup \mathcal{I}_{2}$. \\
	$\mathcal{I}_{1}$ & : &The set of mandatory surgeries (or patients) with the latest surgery days within the planning horizon.  \\
	$\mathcal{I}_{2}$ & : &  The set of optional surgeries (or patients) with the latest surgery days out of the planning horizon.\\
	$\mathcal{S}$ & :& The set of specialties.  \\
	$\mathcal{D}$ & : &  The set of days in the planning horizon.\\
\end{tabular}

	\hspace{-5pt} \begin{tabular}{p{12pt}p{5pt}p{15.45cm}} 	
	$\mathcal{D}_{i}$ & : &  The set of days on which surgery $i$ can be operated considering its time window.\\
	$\mathcal{R}$ & : &  The set of operating rooms.\\
	$\mathcal{R}_{i}$ & : &  The set of operating rooms in which surgery $i$ can be operated in the case that operating rooms are not identical and special equipment is 
	available in some operating rooms.\\
	$\mathcal{H}$ & : &  The set of downstream units. Here, we consider $\mathcal{H}=\{1,2\}$ where 1 and 2 refer to ICU and wards, respectively. It is a general setting that can consider cases with more than two post-surgical recovery areas.\\
\end{tabular}

\vspace{5pt}

\noindent
Parameters:\vspace{6pt}

	\hspace{-5pt} \begin{tabular}{p{30pt}p{5pt}p{14.8cm}} 
	$M_{h}$ & : & The total number of beds in downstream $h \in \mathcal{H}$. \\
	$\alpha _{h}^{shared}$ & : & The percentage of beds in downstream $h$ that are shared between specialties.  \\
	$\alpha _{i}^{waiting}$ & : &   The daily waiting cost of patient $i$.\\
	$e_{i}$ & :& The earliest day in the time window of surgery $i$.  \\
	$s_{i}$ & : &  The specialty corresponding to surgery $i$.\\
\end{tabular}

\hspace{-5pt} \begin{tabular}{p{30pt}p{5pt}p{14.8cm}} 		
	$c_{id}^{waiting}$ & : &  The waiting cost of surgery $i$ if it is performed on day $d$. It is pre-computed by $c_{id}^{waiting}=\alpha _{i}^{waiting}(d-e_{i})$.\\
	$c_{i}^{postpone}$ & : &  The cost of postponing optional surgery $i$ to the next planning horizon.\\
	$c^{OR}$ & : &  The fixed cost of opening an operating room.\\
	$A_{s}^{\min }$ & : &  The minimum number of operating rooms that must be allocated to specialty $s$ in the planning horizon.\\
	$A_{s}^{\max }$ & : &  The maximum number of operating rooms that can be allocated to specialty $s$ in the planning horizon.\\
\end{tabular}

\vspace{5pt}
	
\noindent
Variables:

\vspace{5pt}

\hspace{-5pt} \begin{tabular}{p{12pt}p{5pt}p{15.8cm}} 
	$x_{ird}$ & : & 1 if we allocate surgery $i$ to operating room $r$ on day $d$; 0 otherwise. \\
	$x_{i}^\prime$ & : & 1 if we postpone surgery $i$ to the next planning horizon; 0 otherwise.  \\
\end{tabular}

\hspace{-5pt} \begin{tabular}{p{12pt}p{5pt}p{15.8cm}} 	
	$y_{rd}$ & : &  1 if we open operating room $r$ on day $d$; 0 otherwise.\\
	$z_{srd}$ & : &  1 if specialty $s$ is allocated to operating room $r$ on day $d$; 0 otherwise.\\
	$u_{sh}$ & : &  The number of beds in downstream $h$ allocated to specialty $s$.\\
\end{tabular}

\vspace{15pt}
\noindent
Based on the given notation, we introduce the first-stage model as follows:

\vspace{-8pt}

	\begin{subequations}\label{Model1}
	\renewcommand{\theequation}{\arabic{equation}}
	\begin{align}
	&\hspace{0pt} \min_{x,x^\prime,y,z,u} \sum_{i \in \mathcal{I}}{\sum_{d \in \mathcal{D}_{i}}{\sum_{r \in \mathcal{R}_{i}}{c_{id}^{waiting}x_{ird}}}}+\sum_{i \in \mathcal{I}_{2}}{c_{i}^{postpone}x_{i}^\prime}+\sum_{r \in \mathcal{R}}{\sum_{d \in \mathcal{D}}{c^{OR}y_{rd}}}+Q(x, 
	x^\prime, y, z,u)	&& \label{B1} \\[3pt]
	& \hspace{10pt}  \small \text{Subject to:} & \hspace{0 pt}  \nonumber \\[-5pt]
	&\hspace{120pt} \sum_{d \in \mathcal{D}_{i}}{\sum_{r \in \mathcal{R}_{i}}{x_{ird}}}=1  && \hspace{-155pt} i \in \mathcal{I}_{1} \label{B2} \\[3pt]	
	&\hspace{107pt} \sum_{d \in \mathcal{D}_{i}}{\sum_{r \in \mathcal{R}_{i}}{x_{ird}}}+ x_{i}^\prime=1 && \hspace{-155pt} i \in \mathcal{I}_{2} \label{B3} \\[3pt]
	&\hspace{127pt} \sum_{s \in \mathcal{S}}{z_{srd}}=y_{rd} && \hspace{-155pt} r \in \mathcal{R}, d \in \mathcal{D}  \label{B4} \\[3pt]
	&\hspace{134pt} x_{ird} \le z_{s_{i}rd} && \hspace{-155pt} i \in \mathcal{I}, r \in \mathcal{R}_{i}, d \in \mathcal{D}_{i} \label{B5} \\[3pt]	
	&\hspace{98pt} A_{s}^{\min } \le \sum_{d \in \mathcal{D}}{\sum_{r \in \mathcal{R}}{z_{srd}}}\le A_{s}^{\max } && \hspace{-155pt} s \in \mathcal{S} \label{B6} \\[3pt]
	&\hspace{97pt} \sum_{s \in \mathcal{S}}{u_{sh}} \le \lceil(1-\alpha _{h}^{shared})M_{h}\rceil && \hspace{-155pt} h \in \mathcal{H} \label{B7}\\[3pt]
	&\hspace{129pt} x_{ird}  \in \{0,1\}  && \hspace{-155pt} i \in \mathcal{I}, r \in \mathcal{R}_{i}, d \in \mathcal{D}_{i} \label{B7}\\[3pt]
	&\hspace{133pt} x_{i}^\prime \in \{0,1\} && \hspace{-155pt} i \in \mathcal{I}_{2} \label{B7}\\[3pt]
	&\hspace{131pt} y_{rd}  \in \{0,1\}  && \hspace{-155pt} r \in \mathcal{R}, d \in \mathcal{D} \label{B7}\\[3pt]
	&\hspace{129pt} z_{srd}  \in  \{0,1\} && \hspace{-155pt} s  \in \mathcal{S},r \in \mathcal{R}, d \in \mathcal{D} \label{B7}\\[3pt]
	&\hspace{120pt} u_{sh} \ge 0, integer &&  \hspace{-155pt} s \in \mathcal{S}, h \in \mathcal{H} \label{B7}
	\end{align}
	\end{subequations}

Objective function (1) consists of four components. The first three components calculate patients' waiting cost, postponement cost, and the opening cost of opening operating rooms, respectively. In the last component, we have $Q(x, x^\prime, y, z,u)$ that represents the expected second-stage cost. Constraint (2) indicates that each mandatory surgery must be allocated to one operating room on a single day in the planning horizon. Constraint (3) implies that optional surgeries are either allocated to a surgical block or postponed to the next planning horizon. Constraint (4) states if an operating room is opened on a day, it must be assigned to exactly one specific specialty. Constraint (5) guarantees that surgery $i$ can be operated in operating room $r$ on day $d$ only if this block is assigned to the specialty of surgery $i$ denoted by $s_{i}$. Constraint (6) restricts the number of operating rooms that each specialty can have in the planning horizon. Constraint (7) indicates that, in each downstream, the number of non-shared beds allocated to different specialties cannot be more than the total number of available non-shared beds. In this constraint, $\lceil (1-\alpha _{h}^{shared})M_{h} \rceil$ indicates the smallest integer value larger than or equal to $(1-\alpha _{h}^{shared})M_{h}$.

It is worth noting that, in the first-stage model, determining non-shared beds for each specialty through the variable $u_{sh}$ makes sense only if the planning horizon is long enough. In practice, the planning horizon must be at least one month for the model to be applicable. In this work, we have considered the possibility bed allocation to specialty within a planning horizon to explore the benefits of such a policy.

\subsection{Second-stage model}

In the following, we present the sets, parameters, and variables used in 
the second-stage model.

\vspace{0pt}
\noindent
Sets:\vspace{6pt}

\hspace{-5pt} \begin{tabular}{p{12pt}p{5pt}p{15.45cm}} 
	$\Omega$ & : & The set of stochastic scenarios. \\
	$\mathcal{I}_{s}$ & : &The set of surgeries that belong to specialty $s$.  \\
\end{tabular}

\vspace{20pt}

\noindent
Parameters:\vspace{6pt}

\hspace{-5pt} \begin{tabular}{p{30pt}p{5pt}p{14.8cm}} 
	$c_{h}^{bed}$ & : & The per-unit cost of surge capacity in downstream $h$. \\
	$c^{overtime}$ & : & The per-minute cost of overtime in an operating room.  \\
	$A$ & : &  The total regular available time in each operating room.\\
	$O^{max}$ & :& The maximum allowed overtime in each operating room on each day.  \\
	$t_{i\omega }$ & : &  The duration of surgery $i$ in scenario $\omega $.\\
	$l_{ih\omega }$ & : & The length of stay for patient $i$ in downstream $h$ in 
	scenario $\omega $.\\
	$b_{shd\omega}$ & : & The number of beds that are occupied by existing patients of 
	specialty $s$ in downstream $h$ on day $d$ in scenario $\omega $, from the previous planning horizon.  \\
	$\xi (\omega )$ & : &  The vector of uncertain parameters including surgical durations and lengths of stay in scenario $\omega $.\\
\end{tabular}

\vspace{13pt}
\noindent
Variables:\vspace{6pt}

\hspace{-5pt} \begin{tabular}{p{12pt}p{5pt}p{15.8cm}} 
	$o_{rd\omega }$ & : & The overtime of operating room $r$ on day $d$ in scenario $
	\omega $. \\
	$q_{shd\omega }$ & : & The number of shared beds that are occupied by patients of 
	specialty $s$ in downstream $h$ on day $d$ in scenario $\omega $.  \\
	$v_{shd\omega }$ & : &  The number of patients belonging to specialty $s$ that use the 
	surge capacity in downstream $h$ on day $d$ in scenario $\omega $.\\
\end{tabular}

\vspace{15pt}
\noindent
We formulate the second-stage model as follows:

\vspace{-15pt}

	\begin{subequations}\label{Model1}
	\renewcommand{\theequation}{\arabic{equation}}
	\setcounter{equation}{12}
	\begin{align}
	&\hspace{60pt} Q(x, x^\prime, y, z,u,\omega )= \min 
	_{q,v,o}\sum_{s \in \mathcal{S}}{\sum_{h \in \mathcal{H}}{\sum_{d \in \mathcal{D}}{c_{h}^{bed}v_{shd\omega 
	}}}}+\sum_{d \in \mathcal{D}}{\sum_{r \in \mathcal{R}}{c^{overtime}o_{rd\omega }}}	&& \label{B1} \\[3pt]
	& \hspace{10pt}  \small \text{Subject to:} & \hspace{0 pt}  \nonumber \\[-5pt]
	&\hspace{-5pt} b_{shd\omega} + \sum_{i \in \mathcal{I}_{s}}{\sum_{r \in \mathcal{R}_{i}}{\sum\limits_{\substack{d^\prime \in \mathcal{D}_{i}: d^\prime+\sum_{h^\prime=1}^{h-1}{l_{ih^\prime\omega }} \le d \, \& \, \\ d^\prime+\sum_{h^\prime=1}^{h}{l_{ih^\prime\omega }}> d}}}{x_{id^\prime r}}} \le u_{sh}+q_{shd\omega }+v_{shd\omega }  && \hspace{-80pt} d \in \mathcal{D}, h \in \mathcal{H}, s  \in \mathcal{S},\omega  \in \Omega \label{B2} \\[3pt]	
	&\hspace{90pt} \sum_{s \in \mathcal{S}}{q_{shd\omega }} \le \lfloor\alpha _{h}^{shared}M_{h}\rfloor && \hspace{-80pt} d \in \mathcal{D}, h \in \mathcal{H}, \omega  \in \Omega \label{B3} \\[3pt]
	&\hspace{50pt} \sum_{i \in \mathcal{I}:d \in \mathcal{D}_{i} and r \in \mathcal{R}_{i}}{t_{i\omega }x_{idr}} \le A+ o_{rd\omega 
	} && \hspace{-80pt} d \in \mathcal{D}, r \in \mathcal{R},\omega  \in \Omega  \label{B4} \\[3pt]
	&\hspace{100pt} 0\le o_{rd\omega } \le O^{max} && \hspace{-80pt} d \in \mathcal{D}, r \in \mathcal{R}, \omega  \in \Omega \label{B5} \\[3pt]
	&\hspace{115pt} q_{shd\omega } \ge 0 && \hspace{-80pt} d \in \mathcal{D}, h \in \mathcal{H},s  \in \mathcal{S}, \omega  \in \Omega \label{B6} \\[3pt]
	&\hspace{115pt} v_{shd\omega } \ge 0 && \hspace{-80pt} d \in \mathcal{D}, h \in \mathcal{H},s  \in \mathcal{S}, \omega  \in \Omega \label{B7}
	\end{align}
	\end{subequations}

Objective function (13) minimize the total second-stage cost including the surge capacity cost in different downstream areas and the overtime cost in operating rooms. The second-stage cost $Q(x, x^\prime , y, z,u)$ in objective function (1) is calculated by $E_{\omega  \in \Omega }\lbrack Q(x, x^\prime , y, z,u,\xi (\omega ))\rbrack $ where $E_{\omega  \in \Omega }\lbrack .\rbrack $ computes the expected value over scenarios $\omega  \in \Omega $. Constraint (14) implies the restriction on the number of available beds 
in different downstream areas. In this constraint, ${{\sum\limits_{\substack{d^\prime  \in \mathcal{D}_{i}: d^\prime +\sum_{h^\prime =1}^{h-1}{l_{ih^\prime \omega }} \le d \, \& \,  d^\prime +\sum_{h^\prime =1}^{h}{l_{ih^\prime \omega }}> d}}}{x_{id^\prime r}}}$ is equal to 1 if patient $i$, that is operated in operating room $r$, is in downstream $h$ on day $d$. This is because $d^\prime +\sum_{h^\prime =1}^{h-1}{l_{ih^\prime \omega }} \le d$ indicates that the patient has left the previous downstream $(h-1)$ not later than day $d$ and $d^\prime +\sum_{h^\prime =1}^{h}{l_{ih^\prime \omega }} > d$ shows that he/she will leave downstream $h$ after day $d$. Therefore, the left-hand side of Constraint (14) computes the total number of existing and new patients, associated with specialty $s$ , that are in downstream $h$ on day $d$ in scenario $\omega $. Whenever the left-hand side of this constraint is larger than $u_{sh}$, the model prefers to compensate for the shortage by giving positive values to $q_{shd\omega }$ first and then to $v_{shd\omega }$. This is because the shared ICU beds are available free of cost, while the surge beds, denoted by $v_{shd\omega }$, are penalized in the Objective function (13). We note that, in Constraint (14), if a patient does not need to stay in downstream $h$, it will not be counted as one of the patients for that downstream because in this case the condition $d^\prime +\sum_{h^\prime =1}^{h-1}{l_{ih^\prime \omega }} \le d \,\,\, \& \,\,\,  d^\prime +\sum_{h^\prime =1}^{h}{l_{ih^\prime \omega }}> d$ is not satisfied as $l_{ih\omega}$ is equal to 0.

Constraint (15) shows that the total number of shared beds allocated in each downstream cannot be more than the total number of beds available for sharing. In this constraint, $\lfloor\alpha _{h}^{shared}M_{h}\rfloor$ indicates the largest integer value smaller than or equal to $\alpha _{h}^{shared}M_{h}$. Using the round-up and round-down notations in Constraints (7) and (15), respectively, ensures that all available beds are properly accounted for. Constraint (16) declares that a limited regular time $A$ is available in each operating room on each day. If the total surgical time on the left-hand side of Constraint (16) is more than the regular available time, then the overtime is considered by giving a positive value to $o_{rd\omega }$. Furthermore, Constraint (17) confines the allowed overtime. Constraints (18)-(19) declare the non-negativity of $q_{shd\omega }$ and $v_{shd\omega }$. As we will discuss later in Section 4.2., there is always an optimal feasible solution in which these two variables take integer values.

Although considering surge capacity is beneficial, it should be noted that there are related challenges, especially when specialized nursing skills are required for patient care. In units where patient care depends on nurses with unit-specific training, staffing availability may constrain the use of surge beds. In addition, activating surge capacity may strain supporting resources (e.g., equipment, diagnostics, or coordination across units), potentially affecting care delivery if not carefully managed. Nonetheless, surge capacity remains a practical mechanism for absorbing short-term demand increases, particularly when supported by appropriate staffing and operational adjustments.

\section{Solution algorithm}

One of the main challenges in using stochastic programming models is that solving them becomes exponentially more difficult as the number of scenarios increases. Therefore, in this research, we develop a sample average approximation algorithm (SAA). We also proposed a specialized algorithm to solve the second-stage model for any given first-stage solution quickly.

\subsection{Sample average approximation}
Sample average approximation algorithm method uses Monte Carlo simulation to estimate the expected value of objective function based on a number of random independent identically distributed (i.i.d.) samples. It generates these samples iteratively and solves the extensive form of the stochastic programming model for a limited number of scenarios in each iteration separately. Then, the lower and upper bounds of the objective function are estimated using the outputs of all iterations. Here, we use $|N|$ to denote the number of randomly generated scenarios where N = ${\omega _{1}, \ldots , \omega _{|N|}}$ is the set of scenarios. In SAA, we approximate the second-stage objective value

\vspace{-10pt}

	\begin{subequations}\label{Model1}
	\renewcommand{\theequation}{\arabic{equation}}
	\setcounter{equation}{19}
	\begin{align}
	&\hspace{25pt} E_{\omega \in \Omega}\bigg[ 
	\sum_{s \in \mathcal{S}}{\sum_{h \in \mathcal{H}}{\sum_{d \in \mathcal{D}}{c_{h}^{bed}v_{shdw}}}}+\sum_{d \in \mathcal{D}}{\sum_{r \in \mathcal{R}}{c^{overtime}o_{rdw}}}\bigg] 
		&& \label{B1} \\[3pt]
	& \hspace{27pt}  \small \text{Subject to (14)-(19).} & \hspace{0 pt}  \nonumber
	\end{align}
\end{subequations}

\vspace{-10pt}

\noindent
by

\vspace{-35pt}

	\begin{subequations}\label{Model1}
	\renewcommand{\theequation}{\arabic{equation}}
	\setcounter{equation}{20}
	\begin{align}
	&\hspace{25pt} \min \frac{1}{\vert N\vert }\sum_{n \in N}{\bigg[ 
		\sum_{s \in \mathcal{S}}{\sum_{h \in \mathcal{H}}{\sum_{d \in \mathcal{D}}{c_{h}^{bed}v_{shdn}}}}+\sum_{d \in \mathcal{D}}{\sum_{r \in \mathcal{R}}{c^{overtime}o_{rdn}}}\bigg] 
	}
	&& \label{B1} \\[3pt]
	& \hspace{29pt}  \small \text{Subject to (14)-(19).} & \hspace{0 pt}  \nonumber 
	\end{align}
\end{subequations}

\noindent
In (21), constraints (14)-(19) are repeated for $n  \in N$ instead of $\omega  \in  \Omega $.

The privilege of using SAA is that by increasing the number of samples, the model asymptotically obtains the optimal solution and objective value (\citealp{kleywegt2002sample}). The only issue here is that using a massive number of samples works against the initial intention of using SAA since it becomes a time-consuming procedure. To avoid that, we solve the SAA problem $|M|$ times for a reasonable number of samples $|N|$. Then, we calculate the average and variance of the lower and upper bound using the objective values obtained in all iterations. The procedure of SAA is as follows:

\begin{itemize}
	\item We solve Model (21) $\vert M\vert $ times independently. In each 
	iteration, we consider $|N|$ samples of scenarios and save the 
	obtained first-stage solution ${x}_{N}^{m}$ and the objective value $
	{f}_{N}^{m}$ for each iteration $m \in \{1,\ldots ,\vert M\vert \}$. We 
	refer to the problem that we solve in this iteration as the lower bound 
	problem.
	\item We calculate the average and variance for the lower bound of the 
	objective value using the following formulas:
\end{itemize}

\vspace{-30pt}

	\begin{subequations}\label{Model1}
	\renewcommand{\theequation}{\arabic{equation}}
	\setcounter{equation}{21}
	\begin{align}
	&\hspace{55pt} \overline{LB}= \frac{1}{|M|} \sum_{m=1}^{|M|}{{f}_{N}^{m}}, 
	&& \label{B1} \\
	&\hspace{55pt} \sigma _{LB}^{2}= \frac{1}{|M|(\vert M\vert -1)} 
	\sum_{m=1}^{|M|}{({f}_{N}^{m}- \overline{LB})}.  && \label{B2} 
	\end{align}
\end{subequations}

\begin{itemize}
	\item Then, we solve Model (21) again for $|M|$ iterations for the corresponding fixed first-stage solutions ${x}_{N}^{m}$, $m \in \lbrace 1,\ldots ,\vert M\vert \rbrace $ obtained in the previous step. We refer to the problems solved in this step as upper bound problems. In each iteration, we consider $|P|$ samples of scenarios and save ${f}_{p}({x}_{N}^{m})$ that denotes the objective value of the upper bound problem in iteration $m$ for sample $p$. 
	\item Next, the average and the variance for the upper bound of the objective value are calculated by:
\end{itemize}

\vspace{-15pt}

	\begin{subequations}\label{Model1}
	\renewcommand{\theequation}{\arabic{equation}}
	\setcounter{equation}{23}
	\begin{align}
	&\hspace{55pt} \overline{UB}= \frac{1}{|P|} \sum_{p=1}^{|P|}{{f}_{p}({x}_{N}^{m})}, 
	&& \label{B1} \\
	&\hspace{55pt} \sigma _{UB}^{2}= \frac{1}{|P|(\vert P\vert -1)} 
	\sum_{p=1}^{|P|}{({f}_{p}({x}_{N}^{m})- \overline{UB})}.  && \label{B2} 
	\end{align}
\end{subequations}

\begin{itemize}
	\item We consider the minimum of average upper bound values $\overline{UB}$ obtained by (24) over different iterations as the best upper bound value and denote it by $\overline{UB}_{Best}$. Then we use $\overline{LB}$ and $\overline{UB}_{Best}$ to calculate the optimality gap by
\end{itemize}

\vspace{-25pt}

	\begin{subequations}\label{Model1}
	\renewcommand{\theequation}{\arabic{equation}}
	\setcounter{equation}{25}
	\begin{align}
	&\hspace{55pt} Gap= 100\frac{(\overline{UB}_{Best}- \overline{LB})}{\overline{LB}}.
	&& \label{B1} 
	\end{align}
\end{subequations}

\vspace{5pt}

Additionally, in the implementation of the lower bound problem, we include the following constraint to prevent infeasibility in the upper bound problem arising from limited operating room overtime availability:

\begin{center}
	$\sum_{i \in \mathcal{I}:d \in \mathcal{D}_{i} and r \in \mathcal{R}_{i}}{t_{i}^{max} x_{idr}} \le A+ O^{max} \hspace{20 pt} d \in \mathcal{D}, r \in \mathcal{R}$
\end{center}

\vspace{5pt}

In the above constraint, $t_{i}^{\max}$ denotes the maximum possible duration of surgery $i$.

\subsection{Specialized algorithm}
In this section, we propose a specialized algorithm that optimally solves the second-stage model (13)-(19) for any given first-stage solution, without any need to solve a linear programming model after fixing the integer variables (e.g., see \citealp{hashemi2012effective}). This specialized algorithm is very quick and therefore can solve the upper bound problem of the sample average approximation for all scenarios in a very short time. As we will see in Section 5, thanks to this algorithm, we can evaluate a very large number of scenarios in the upper bound problems and obtain reliable statistical upper bounds with small standard deviations. We reiterate that the proposed specialized algorithm is used only to solve the upper-bound problem and not the lower-bound problem.

We know that, on each day $d$, at most $u_{sh}$ patients can be allocated to non-shared beds of specialty s. Therefore, in in each scenario $\omega$ and each downstream level $h$, the total number of patients that cannot be allocated to the non-shared beds is equal to \\ 

\begin{center}
	$max \Biggl\{b_{shd\omega} + \sum_{i \in \mathcal{I}_{s}}{\sum_{r \in \mathcal{R}_{i}}{\sum\limits_{\substack{d^\prime \in \mathcal{D}_{i}: d^\prime+\sum_{h^\prime=1}^{h-1}{l_{ih^\prime\omega }} \le d \, \& \, \\ d^\prime+\sum_{h^\prime=1}^{h}{l_{ih^\prime\omega }}> d}}}{x_{id^\prime r}}} - u_{sh},0\Biggr\}$.
\end{center}

These patients must be allocated to either the shared beds or must be served by the surge capacity. Since using shared beds is free of cost, the model always first tries to allocate the patients to shared beds. If they are not enough, the remaining patients must be served by the surge capacity. In the allocation of shared beds to specialties, the proposed stochastic programming model is indifferent to different specialties because the cost of using surge capacity is the same for all of them. Therefore, we can randomly allocate the patients not served by non-shared beds, using the following procedure: First, we produce a random order of specialties denoted by $s_{1},s_{2},…,s_{|S|}$. Then, we use the following formula to assign the shared beds to specialties:

\vspace{-10pt}

\begin{equation*}
\setcounter{equation}{27}
\resizebox{\textwidth}{!}{$\displaystyle % restart math mode!
	q_{s_{k}hd\omega} =
	\begin{dcases}
	0 & \text{if $\alpha _{h}^{shared}M_{h} - \sum_{k^\prime=1}^{k-1}{q_{s_{k\prime}hd\omega }} \le 0$}\\
	min \Bigggggl\{ max \Biggggl\{  b_{shd\omega} + \sum_{i \in \mathcal{I}_{s_{k}}}{\sum_{r \in \mathcal{R}_{i}}{\sum\limits_{\substack{d^\prime \in \mathcal{D}_{i}: d^\prime+\sum_{h^\prime=1}^{h-1}{l_{ih^\prime\omega }} \le d \, \& \, \\ d^\prime+\sum_{h^\prime=1}^{h}{l_{ih^\prime\omega }}> d}}}{x_{id^\prime r}}} - u_{sh},0 \Biggggr\},  \alpha _{h}^{shared}M_{h} - \sum_{k^\prime=1}^{k-1}{q_{s_{k\prime}hd\omega }}    \Bigggggr\} & \text{if $\alpha _{h}^{shared}M_{h} - \sum_{k^\prime=1}^{k-1}{q_{s_{k\prime}hd\omega }} > 0$}
	\end{dcases}
	$} % end math mode and close the box to be resized
\end{equation*}

\vspace{3pt}

\begin{flushright}
	$k \in \{1,...,|\mathcal{S}|\}, h \in \mathcal{H}, d \in \mathcal{D}, \omega \in \Omega$
\end{flushright}

In the above formula, if $\alpha _{h}^{shared}M_{h} - \sum_{k^\prime=1}^{k-1}{q_{s_{k\prime}hd\omega }} \le 0$ happens, it means the previous specialties in the random list have already used all shared beds, therefore the remaining patients of specialty $s_k$ must be passed to the surge capacity. Otherwise, we will allocate some beds to specialty. In this case, the number of shared beds allocated to specialty $s_k$ cannot be more than the number of remaining patients of this specialty that could not be allocated to non-shared beds, and the number of remaining shared beds. That is why $q_{s_{k}hd\omega}$ is computed in the format of min{.,.} in the second relation.

The above formula computes the number of shared beds that must be allocated to each specialty for different downstream and days in the planning horizon separately. After computing the values of $q_{shd\omega}$, we can simply determine the values of $v_{shd\omega}$ by the following relation:

\vspace{5pt}

$v_{s_{k}hd\omega} = max \Biggl\{b_{shd\omega} + \sum_{i \in \mathcal{I}_{s_{k}}}{\sum_{r \in \mathcal{R}_{i}}{\sum\limits_{\substack{d^\prime \in \mathcal{D}_{i}: d^\prime+\sum_{h^\prime=1}^{h-1}{l_{ih^\prime\omega }} \le d \, \& \, \\ d^\prime+\sum_{h^\prime=1}^{h}{l_{ih^\prime\omega }}> d}}}{x_{id^\prime r}}} - u_{s_{k}h}-q_{s_{k}hd\omega},0\Biggr\}$

\begin{flushright}
	$k \in \{1,...,|\mathcal{S}|\}, h \in \mathcal{H}, d \in \mathcal{D}, \omega \in \Omega$
\end{flushright}

%For the sake of fairness and to avoid giving a full priority to a specialty for the whole horizon, in using the above formulas, we produce different random orders of specialties (i.e., $s_1,s_2,…,s_{|S|}$) for each scenario, each day, and each downstream level.
It is noteworthy that the proposed specialized algorithm shows the existence of an optimal solution in the second-stage model that satisfies the integrality requirement of the variables $q_{shd\omega}$ and $v_{shd\omega}$, as needed in practice.

\vspace{5pt}

\section{Computational results}

We have separated this section into five different subsections. In subsection 5.1, we explain the generation of instances. In subsection 5.2, we computationally evaluate the improvement in the objective value resulting from sharing beds in downstream units using our proposed stochastic programming model. In subsection 5.3, we tune the main parameters of our sample average approximation. Then, in subsection 5.4, we evaluate the performance of our model within the framework of sample average approximation. Finally, in subsection 5.5, we perform an extensive sensitivity analysis for cost parameters to measure the efficiency of our model in different situations.

We used IBM ILOG CPLEX Optimization to solve the integer programming model we implemented our code in Visual Studio V12.8 in C++. The experiments were run on a computer with two AMD Rome 7502 processors, 2.50 Ghz, and a total of 64 cores. We used a single core to run each instance. All source code and generated instances are available at \url{https://bit.ly/4bL572C}.

\subsection{Instances}
To generate instances, we set the number of weeks in the planning horizon to \{2,3,4\}. In our instances, four operating rooms are available on each day. We considered eight hours as the daily regular available time of each operating room. Also, the maximum allowed overtime for each operating room is set to three hours. Furthermore, we have considered ICU and wards as two consecutive downstream levels. We set the number of available ICU and ward beds to 35 and 65, respectively. For instances with two, three, and four weeks, we set the number of patients to \{120,180,240\}, respectively. The earliest days of time windows for surgeries are randomly generated within the planning horizon. Also, the length of the time window for each patient varies from one to seven days. Each patient is randomly associated with one of the seven specialties listed in Table 1. This table also reports the average and standard deviation of surgical times, as well as the average and standard deviation of lengths of stay in wards and the ICU for each specialty (\citealp{costa2017assessment}). Additionally, to set the total number of available beds (i.e.,$M_h$), we first compute the sum of the average lengths of stay for all patients, divide this value by the number of days in the planning horizon, and then multiply the result by 0.8. This ensures that the number of beds is lower than the total expected demand, making the instances both computationally challenging and representative of realistically limited resources. We have provided more details about the generation of cost and uncertain parameters of the instances in Appendix 1.

To generate our instance set, we considered that the number of specialties in each instance belongs to \{1,2,3,4,5,6,7\}. Considering all possibilities for the number of weeks in the planning horizon and the number of specialties, we have 21 combinations. For each combination, we generated 5 instances for a total of 105 instances.

	\begin{table}[bt]
	%\caption{\hspace{-5pt} The list of specialties, the average surgical time, the average and standard deviation of lengths of stay (LOS).}	
	\label{Table2}
	\resizebox{\textwidth}{!}{
		\small
		\begin{tabular}{lccccc}			
			\multicolumn{6}{l}{\textbf{Table 1} \,\,\, The list of specialties, the average surgical time, the average and standard deviation of lengths of stay (LOS).}\\
			\hline\noalign{\smallskip}
			Specialty & \specialcell{Average surgical\\ time (min)} & \specialcell{Standard deviation of\\ surgical time (min)}  & \specialcell{Average LOS \\in wards (day)} & \specialcell{Average LOS\\ in ICU (day)} & \specialcell{Standard deviation of LOS \\in both wards and ICU (day)}\\
			\hline\noalign{\smallskip}
			General	&	150.95	&	25.16	&	3.10	&	4.65	&	4.48	\\[2pt]		
			Neurology	&	135.06	&	22.51	&	2.89	&	4.34	&	5.19	\\[2pt]		
			Cardiovascular	&	189.34	&	31.56	&	2.34	&	3.50	&	3.01	\\[2pt]		
			Orthopedic	&	151.95	&	25.33	&	3.08	&	4.61	&	4.51	\\[2pt]		
			Urology	&	94	&	5.22	&	6.27	&	9.402	&	3.68	\\[2pt]		
			Plastic	and	reconstructive	&	157.72	&	10.52	&	15.77	&	6.71	&	4.54	\\[2pt]
			Obstetrics	and	gynecology	&	79.32	&	5.29	&	7.93	&	5.22	&	2.21	\\[2pt]
			\hline\noalign{\smallskip}
		\end{tabular}
	}
\end{table}

\subsection{Analysis of bed-sharing policies}
In this section, we analyze the improvement in the objective value obtained by sharing beds in downstream units. We have considered three levels of bed sharing in our computational analysis. In the first one, referred to as the ``No-Sharing'' policy, we suppose that sharing of downstream beds is not allowed among specialties. This means that $\alpha _{h}^{shared}$ is set to 0 for $h=1,2$. In the second setting, we considered a ``Midlevel-Sharing'' policy where 50\% of all available beds are shared among specialties. This setting refers to $\alpha _{h}^{shared}=0.5$ for $h=1,2$. In the last setting, referred to as the ``Full-Sharing'' policy, we assume that all beds are shared among specialties without any limitation. In this case, we have $\alpha _{h}^{shared}=1$ for $h=1,2$.

To obtain the results of these three pooling policies, we solved the extensive form of the proposed stochastic programming model without SAA, so that all policies were evaluated under the same set of 30 scenarios. This avoids discrepancies caused by the random scenario generation and statistical variability inherent in SAA outputs. We have presented the results of the above sharing policies in Table 2. In this Table, we have 18 combinations of instances with different values of the number of weeks and the number of specialties presented under the first two columns. We have ignored instances with one specialty as sharing beds is meaningless in this case. Each row of Table 2 presents the average results of five instances.

In Table 2, we have presented the results of the Midlevel-Sharing (setting 2) and Full-Sharing (setting 3) policies in comparison to the No-Sharing policy (setting 1) separately. The columns of ``Imp. (\%)'' denote the total improvement in the objective value obtained by the corresponding sharing policy (Midlevel-Sharing or Full-Sharing) in comparison to the No-Sharing policy. For instance, under the Column of ``Midlevel-Sharing Policy'', we have $Imp. (\%)=100(Obj_{No sharing}-Obj_{Midlevel})/Obj_{No sharing}$, where $Obj_{No sharing}$ and $Obj_{Midlevel}$ denote the objective values of Midlevel-Sharing and No-Sharing policies, respectively. Similarly, Columns ``Overtime cost Imp. (\%)'', ``Surge cost Imp. (\%)'', ``Waiting cost Imp. (\%)'', ``Postponement cost Imp. (\%)'', and ``OR cost Imp. (\%)'' indicate the contributions of each type of cost in the total improvement. This means that the sum of the recent five columns is equal to the value of ``Imp. (\%)''.

The results of Table 2 demonstrate that the higher volume of sharing results in more improvement in the objective function. In the largest instances with seven specialties and four weeks, we observe that the values of ``Imp. (\%)'' are 16.97\% and 19.53\% for the Midlevel-Sharing and Full-Sharing policies, respectively. Besides, the average values of ``Imp. (\%)'' in the last row of Table 2 show that the Midlevel-Sharing and Full-Sharing policies lead to 11.29\% and 12.38\% improvement compared to the No-Sharing policy. The results also show that the average improvement from Midlevel-Sharing to Full-Sharing is marginal around 1.09\% and therefore the Midlevel-Sharing policy could be enough to significantly improve the performance of operating room planning. 

\afterpage{

	\begin{landscape}
		$\,$
		\vspace{0pt}
		\setlength\LTcapwidth{\textwidth} % default: 4in (rather less than \textwidth...)
		\setlength\LTleft{-15pt}            % default: \parindent
		\setlength\LTright{0pt}           % default: \fill
		{\setlength{\extrarowheight}{3 pt}%
			\small
			\begin{longtable}{@{\extracolsep{\fill}}*{16}{c}}				
				%\caption{}
				%\label{my-label} \\ [-20 pt]
				\multicolumn{16}{l}{\textbf{Table 2} \,\,\, Comparison of bed-sharing policies.}\\
				\hline\noalign{\smallskip}	
				\multicolumn{2}{c}{\textit{Data Info.}}  & & \multicolumn{6}{c}{\textit{Midlevel-Sharing Policy}}  & & \multicolumn{6}{c}{\textit{Full-Sharing Policy}}\\
				\cline{1-2} \cline{4-9} \cline{11-16}\noalign{\smallskip}
				\textit{\specialcell{No. of \\Weeks}} & \textit{\specialcell{No. of \\Spec.}} & & \textit{\specialcell{Imp. \\(\%)}} & \textit{\specialcell{Overtime \\cost Imp. \\(\%)}} & \textit{\specialcell{Surge \\cost \\Imp.\\ (\%)}} & \textit{\specialcell{Waiting \\cost\\ Imp.\\ (\%)}} & \textit{\specialcell{Postpone-\\ment \\cost Imp. \\(\%)}} &\textit{\specialcell{ OR cost\\ Imp.\\ (\%)}} & & \textit{\specialcell{Imp.\\ (\%)}} & \textit{\specialcell{Overtime\\ cost Imp.\\ (\%)}} & \textit{\specialcell{Surge \\cost \\Imp. \\(\%)}} & \textit{\specialcell{Waiting \\cost \\Imp.\\ (\%)}} & \textit{\specialcell{Postpone-\\ment \\cost Imp. \\(\%)}} & \textit{\specialcell{OR cost \\Imp. \\(\%)}}	\\
				\hline\noalign{\smallskip}
				2 & 2 & & 4.56 & 0.06 & 5.61 & 0.42 & -1.79 & 0.27 & & 4.75 & 0.09 & 
				5.62 & 0.88 & -2.01 & 0.16 \\
				& 3 & & 6.42 & -0.14 & 7.29 & -0.44 & -0.32 & 0.03 & & 6.69 & -0.13 & 
				6.94 & -0.04 & 0.00 & -0.07 \\
				& 4 & & 8.34 & -0.01 & 9.75 & 0.88 & -2.25 & -0.03 & & 8.62 & 0.08 & 
				9.54 & 0.79 & -1.82 & 0.03 \\
				& 5 & & 12.69 & 0.06 & 11.75 & 1.35 & -0.78 & 0.31 & & 13.82 & 0.07 & 
				12.49 & 1.64 & -0.86 & 0.48 \\
				& 6 & & 13.62 & 0.00 & 10.43 & 0.66 & 2.32 & 0.22 & & 15.44 & 0.02 & 
				13.53 & 1.08 & 0.52 & 0.28 \\
				& 7 &   & 15.93 & 0.02 & 13.65 & 0.60 & 1.78 & -0.11 &   & 17.79 & 
				-0.02 & 14.49 & 0.56 & 2.90 & -0.14 \\
				\hline\noalign{\smallskip}
				\multicolumn{3}{c}{\textit{Average}}   & 10.26 & 0.00 & 9.74 & 0.58 & -0.17 & 0.12 &   & 11.18 & 
				0.02 & 10.44 & 0.82 & -0.21 & 0.12 \\
				\hline\noalign{\smallskip}
				3 & 2 & & 4.84 & -0.01 & 5.16 & -0.07 & -0.06 & -0.18 & & 4.90 & 0.00 & 
				7.15 & 0.40 & -2.53 & -0.12 \\
				& 3 & & 8.64 & -0.21 & 8.21 & 0.54 & 0.14 & -0.04 & & 8.91 & -0.18 & 
				6.37 & 0.21 & 2.60 & -0.08 \\
				& 4 & & 12.12 & -0.03 & 11.49 & 0.44 & 0.21 & 0.02 & & 12.62 & -0.06 & 
				13.94 & 0.99 & -2.33 & 0.08 \\
				& 5 & & 13.21 & -0.02 & 10.37 & 0.80 & 1.97 & 0.08 & & 14.35 & 0.00 & 
				11.44 & 1.13 & 1.74 & 0.05 \\
				& 6 & & 14.34 & -0.04 & 11.73 & 0.48 & 2.06 & 0.11 & & 16.08 & -0.07 & 
				12.72 & 0.36 & 2.99 & 0.08 \\
				& 7 &   & 16.85 & -0.02 & 11.03 & -0.22 & 6.05 & 0.00 &   & 19.42 & 
				-0.04 & 13.03 & 0.70 & 5.65 & 0.08 \\
				\hline\noalign{\smallskip}
				\multicolumn{3}{c}{\textit{Average}}  & 11.66 & -0.05 & 9.66 & 0.33 & 1.73 & 0.00 &   & 12.71 & 
				-0.06 & 10.77 & 0.63 & 1.35 & 0.01 \\
				\hline\noalign{\smallskip}
				4 & 2 & & 4.99 & -0.03 & 3.58 & -0.07 & 1.66 & -0.16 & & 5.13 & -0.04 & 
				3.51 & -0.05 & 1.82 & -0.11 \\
				& 3 & & 9.58 & -0.10 & 5.62 & -0.41 & 4.58 & -0.11 & & 9.98 & -0.07 & 
				7.07 & -0.19 & 3.21 & -0.04 \\
				& 4 & & 11.43 & 0.00 & 8.87 & 0.83 & 1.67 & 0.06 & & 12.42 & -0.20 & 
				8.72 & 0.51 & 3.19 & 0.19 \\
				& 5 & & 13.82 & -0.03 & 8.38 & -0.67 & 5.96 & 0.18 & & 15.36 & -0.04 & 
				8.63 & -0.49 & 7.14 & 0.12 \\
				& 6 & & 14.90 & -0.02 & 11.04 & -0.11 & 3.93 & 0.05 & & 16.99 & -0.02 & 
				10.82 & -0.19 & 6.35 & 0.04 \\
				& 7 &   & 16.97 & -0.01 & 11.28 & -0.27 & 6.05 & -0.08 & & 19.53 & 
				-0.09 & 12.48 & -0.36 & 7.59 & -0.08 \\[-3pt]
				\hline\noalign{\smallskip}\\[-16pt]
				\multicolumn{3}{c}{\textit{Average}}  & 11.95 & -0.03 & 8.13 & -0.12 & 3.97 & -0.01 &   & 13.24 & 
				-0.08 & 8.54 & -0.13 & 4.88 & 0.02 \\[-3pt]
				\hline\noalign{\smallskip}
				\multicolumn{3}{c}{\textit{Total Average}}   & 11.29 & -0.03 & 9.18 & 0.26 & 1.84 & 0.03 &   & 
				12.38 & -0.04 & 9.92 & 0.44 & 2.01 & 0.05 \\[-3pt]
				\hline
		\end{longtable}}
	\end{landscape}
}

Moreover, we observe a strictly increasing trend of improvement in each set of weeks as the number of specialties rises. For example, the results of the Midlevel-Sharing policy show that the average of total cost improvement changes from 4.56\%, 4.84\%, and 4.99\% to 15.93\%, 16.85\%, and 16.97\% as the number of specialties increases from 2 to 7 in instances with 2, 3, and 4 weeks, respectively. Besides, it is noteworthy that, for a fixed number of specialties, we have a higher value of improvement in longer planning horizons. 

Finally, the values under columns ``Surge cost Imp. (\%)'' and ``Postponement cost Imp. (\%)'' have the highest role in the total cost improvement. More specifically, in Midlevel- and Full-sharing policies, the average saving in the surging capacity are 9.18\% and 9.92\%, respectively. Also, the average saving obtained in postponement cost is 1.84\% and 2.01\%, respectively. For the remaining sections of this research, we run the computational experiments for the case of the Midlevel-Sharing policy.

\subsection{Parameter tuning of sample average approximation method}

In this section, we carry out some computational experiments to tune the parameters of the sample average approximation algorithm. These parameters include the number of iterations $(|M|),$ the number of samples in the lower bound problem $(|N|)$, and the number of samples in the upper bound problem $(|P|).$ To do so, we perform some computational analysis on one of the instances with three weeks of planning horizon and seven specialties. The reason for choosing one of the instances of this combination is that our primary computational results showed that these instances are not too computationally demanding and their medium size lets them mimic the computational behavior of other larger and smaller instances. Therefore, we first ran the model for this instance for different combinations of $|N| \in \{5,10,20,30,40,50,60,70\}$ and $|M| \in \{5,10,15,20,25,30,40,50\}$, while keeping $|P|=30$.

We compared the combinations based on three main features: total solution time, gap, and relative standard deviation. In Figure 1, we study the trend of the optimality gap for different values of $|N|$ and $|M|$. This figure shows that for almost all the values of $|N|$, the amount of gap decreases as the number of $|M|$ increases. It is worth mentioning that the fluctuations of the gap value considerably fall after the value of $|N|$ crosses 30. 

In Figure 2, we observe the behavior of Relative Standard Deviation (RSD) under different circumstances. RSD is a feature that measures the significance of the deviation in the lower bound with respect to its average value and is calculated by $RSD= SD_{LB}/LB$. In this formula,$ LB$ and $SD_{LB}$ denote the estimated lower bound and its standard deviation. Figure 2 shows that, after crossing $|N|=20$, the value of $RSD$ is decreasing gradually as the number of iterations grows. Finally, in Figure 3, we observe that with growth in the number of iterations, the solution time of the model increases. 

Based on the above analysis, we set $|N| = 30$ and $|M| = 25$, as this parameter configuration yields stable SAA results with a small optimality gap (i.e., less than 0.1\%) and RSD, while maintaining a reasonable solution time (approximately 4,100 seconds).

Next, we have to tune the number of samples in the upper bound problem that is denoted by $|P|$. For obtaining the appropriate value of $|P|$, first we fix $\vert N\vert =30$ and $\vert M\vert =25$ and then run the selected instance for different values of $|P|$. Figure 4 depicts $RSD$ and solution time of the upper bound problem for different values of $|P|$. Based on Figure 4, we choose $\vert P\vert =6000$ as it results in a very low value of RSD with a reasonable solution time.

 \begin{figure}[H]
	\begin{center}
		\hspace{-10pt}\includegraphics[height=2.3in]{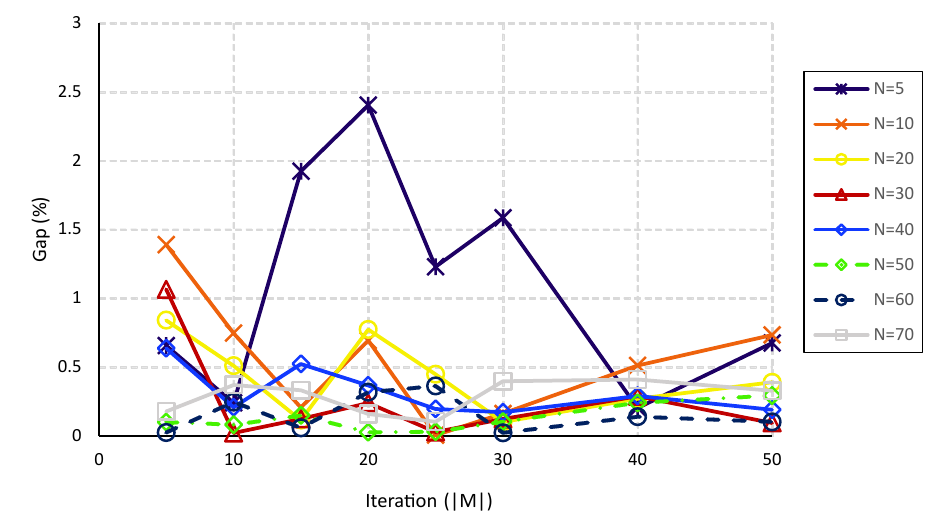}
		\caption{Optimality gap (\%) for different combinations of $|M|$ and $|N|$.} \label{FigureE.C10-1}
	\end{center}
\end{figure}

 \begin{figure}[H]
	\begin{center}
		\hspace{-14pt}\includegraphics[height=2.3in]{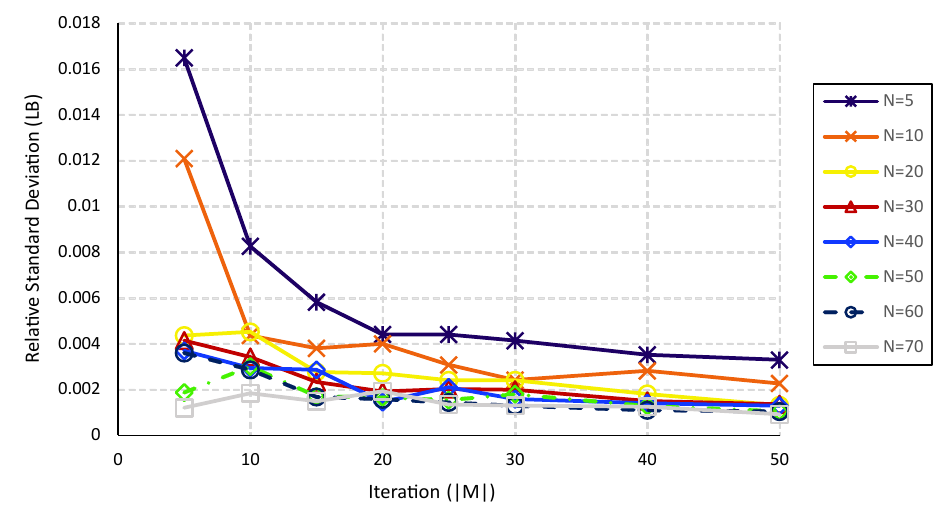}
		\caption{The RSD of LB for different combinations of $|M|$ and $|N|$.} \label{FigureE.C10-1}
	\end{center}
\end{figure}

 \begin{figure}[H]
	\begin{center}
		\hspace{-20pt}\includegraphics[height=2.3in]{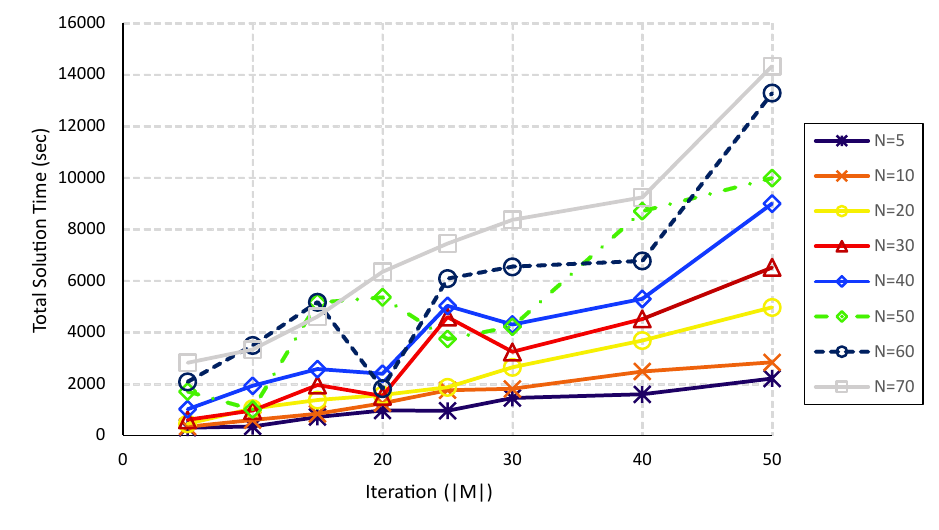}
		\caption{Solution times (sec) for different combinations of $|M|$ and $|N|$.}\label{FigureE.C10-1}
	\end{center}
\end{figure}

 \begin{figure}[t]
	\begin{center}
		\hspace{-33pt}\includegraphics[height=2.3in]{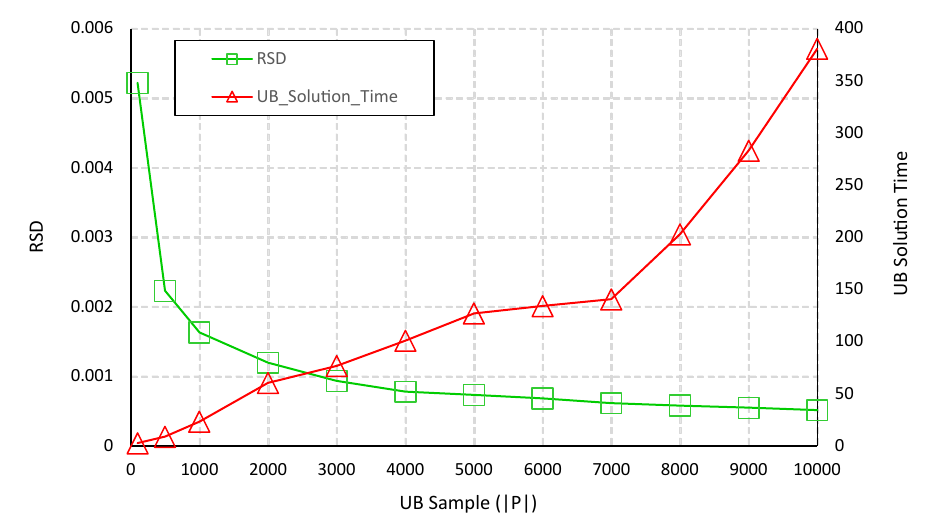}
		\caption{The RSD and solution time of the upper bound problem (sec) for different values of $|P|$.} \label{FigureE.C10-1}
	\end{center}
\end{figure}

\subsection{Performance of SAA algorithm}

In this section, we evaluate the efficiency of the proposed stochastic programming model embedded in the SAA framework using the tuned parameters. In Table 3, we have provided the results of instances with up seven specialties and four weeks of the planning horizon. 

In this table, each row presents the average of five instances for a fixed number of weeks and specialties. Columns ``No. of Weeks'' and ``No. of Spec.'' show the number of weeks in a planning horizon and the number of specialties, respectively. Column ``No. of Ite.'' shows the number of iterations repeated by the SAA algorithm. Column ``Time (Sec)'' indicates the total solution time of the SAA algorithm. Columns ``LB Time (\%)'' and ``UB Time (\%)'' stand for the percentages of total solution time spent to solve the lower bound and upper bound problems, respectively. Columns ``LB'' and ``$SD_{LB}$'' indicate the average and standard deviation of the lower bound. Similarly, columns ``UB'' and ``$SD_{UB}$'' represent the average and standard deviation of the upper bound problem, respectively. The next five columns show the contribution of different cost components in the obtained objective function. Columns ``Overtime cost (\%)'', ``Surge capacity cost (\%)'', ``Waiting cost (\%)'', ``Postponement cost (\%)'', and ``OR cost (\%)'' represent the operating rooms' overtime cost, the total surge cost in downstream units, the patients' waiting cost, the cost of postponing surgeries, and the fixed cost of opening operating rooms, respectively. Under column ``Gap (\%)'', we have provided the optimality gap calculated by formula $100(UB-LB)/LB$. Finally, ``VSS (\%)'' gives the value of stochastic solution that is computed by $100(UB_{EVP}-UB)/UB_{EVP}$. Here, $UB_{EVP}$ is the objective value of the solution obtained by solving the Expected Value Problem (EVP) and then evaluated by the scenarios of the stochastic problem. Expected Value Problem (EVP) is the deterministic version of our stochastic programming model where all random parameters are replaced by their corresponding average values from different scenarios. To obtain the amount of $UB_{EVP}$, we need to run the deterministic EVP for an average scenario. Then, we evaluate the obtained first-stage solution using the same stochastic scenarios. Generally, ``VSS'' determines how better the stochastic programming model works compared to its deterministic EVP in terms of the objective value.

\afterpage{

	\begin{landscape}
		$\,$
		\vspace{0pt}
		\setlength\LTcapwidth{\textwidth} % default: 4in (rather less than \textwidth...)
		\setlength\LTleft{-15pt}            % default: \parindent
		\setlength\LTright{0pt}           % default: \fill
		{\setlength{\extrarowheight}{1.5 pt}%
			\small
			\begin{longtable}{@{\extracolsep{\fill}}*{18}{c}}				
				%\caption{}
				%\label{my-label} \\ [-20 pt]
				\multicolumn{18}{l}{\textbf{Table 3} \,\,\, Computational results of the SAA algorithm.}\\
				\hline\noalign{\smallskip}	
				\multicolumn{2}{c}{\textit{Data Info.}}  & & \multicolumn{15}{c}{\textit{Sample Average Approximation}}\\
				\cline{1-2} \cline{4-18} \noalign{\smallskip}
				\textit{\specialcell{No. of \\Weeks}} & \textit{\specialcell{No. \\ of \\Spec.}} & & \textit{\specialcell{No. of \\Ite.}} & \textit{\specialcell{Time\\ (Sec)}} & \textit{\specialcell{LB Time \\(\%)}} & \textit{\specialcell{UB Time \\(\%)}} & \textit{LB} &\textit{$SD_{LB}$} & \textit{UB} & \textit{$SD_{UB}$} & \textit{\specialcell{Gap \\(\%)}} & \textit{\specialcell{VSS \\(\%)}} & \textit{\specialcell{Overtime \\cost (\%)}} & \textit{\specialcell{Surge\\ capacity \\cost \\(\%)}}& \textit{\specialcell{Waiting\\cost \\(\%)}}& \textit{\specialcell{Postpone-\\ment \\cost \\(\%)}}& \textit{\specialcell{OR \\cost \\(\%)}}	\\
				\hline\noalign{\smallskip}
				2 & 1 & & 25 & 554 & 30.0 & 70.0 & 37319 & 97 & 37406 & 33 & 0.24 & 1.71 
				& 0.37 & 57.81 & 6.38 & 12.47 & 22.96 \\
				& 2 & & 25 & 852 & 62.3 & 37.7 & 38494 & 97 & 38555 & 34 & 0.16 & 3.71 & 
				0.08 & 57.96 & 9.00 & 9.76 & 23.20 \\
				& 3 & & 25 & 34038 & 98.6 & 1.4 & 39959 & 83 & 40113 & 32 & 0.39 & 5.16 
				& 0.72 & 52.37 & 10.80 & 11.25 & 24.85 \\
				& 4 & & 24.8 & 32498 & 98.5 & 1.5 & 39530 & 86 & 39603 & 31 & 0.18 & 
				7.51 & 0.43 & 51.37 & 12.84 & 10.64 & 24.72 \\
				& 5 & & 25 & 46162 & 99.0 & 1.0 & 41826 & 92 & 41884 & 31 & 0.14 & 8.50 
				& 0.35 & 55.11 & 13.43 & 8.26 & 22.84 \\
				& 6 & & 23.2 & 44572 & 99.1 & 0.9 & 43320 & 85 & 43475 & 31 & 0.37 & 
				9.90 & 0.26 & 52.52 & 14.53 & 9.34 & 23.34 \\
				& 7 & & 25 & 17805 & 98.4 & 1.6 & 40956 & 75 & 41076 & 29 & 0.30 & 10.02 
				& 0.16 & 51.51 & 16.19 & 7.33 & 24.81 \\
				\hline\noalign{\smallskip}
				\multicolumn{3}{c}{\textit{Average}} & 25 & 25212 & 84 & 16 & 40200 & 88 & 40302 & 31 & 0.25 & 
				6.64 & 0.34 & 54.09 & 11.88 & 9.86 & 23.82 \\
				\hline\noalign{\smallskip}
				3 & 1 & & 25 & 2380 & 23.7 & 76.3 & 61728 & 153 & 61920 & 47 & 0.31 & 
				1.70 & 0.33 & 55.09 & 5.51 & 18.95 & 20.12 \\
				& 2 & & 25 & 2253 & 49.8 & 50.2 & 67989 & 146 & 68076 & 49 & 0.13 & 3.78 
				& 0.05 & 58.10 & 6.91 & 15.91 & 19.02 \\
				& 3 & & 21.6 & 54109 & 98.1 & 1.9 & 68451 & 127 & 68692 & 45 & 0.36 & 
				10.86 & 0.61 & 55.05 & 8.69 & 15.11 & 20.54 \\
				& 4 & & 23.6 & 61443 & 97.6 & 2.4 & 69761 & 129 & 70013 & 46 & 0.36 & 
				11.30 & 0.51 & 53.77 & 8.32 & 16.55 & 20.85 \\
				& 5 & & 25 & 19033 & 95.1 & 4.9 & 68105 & 115 & 68292 & 43 & 0.27 & 
				11.72 & 0.25 & 54.90 & 9.75 & 13.66 & 21.44 \\
				& 6 & & 25 & 37081 & 97.4 & 2.6 & 71893 & 124 & 72125 & 44 & 0.32 & 
				19.48 & 0.25 & 52.03 & 10.94 & 15.85 & 20.92 \\
				& 7 & & 25 & 17567 & 96.6 & 3.4 & 69732 & 116 & 69810 & 42 & 0.11 & 
				25.74 & 0.18 & 53.54 & 12.49 & 12.50 & 21.29 \\
				\hline\noalign{\smallskip}
				\multicolumn{3}{c}{\textit{Average}} & 24 & 27695 & 80 & 20 & 68237 & 130 & 68418 & 45 & 0.27 & 
				12.08 & 0.31 & 54.64 & 8.94 & 15.50 & 20.60 \\
				\hline\noalign{\smallskip}
				4 & 1 & & 25 & 2625 & 20.4 & 79.6 & 88940 & 155 & 89301 & 57 & 0.41 & 
				1.58 & 0.28 & 50.51 & 4.18 & 27.29 & 17.74 \\
				& 2 & & 25 & 2567 & 43.0 & 57.0 & 98904 & 187 & 99151 & 62 & 0.25 & 5.50 
				& 0.05 & 55.91 & 4.69 & 22.15 & 17.19 \\
				& 3 & & 25 & 36284 & 93.6 & 6.4 & 93645 & 138 & 93815 & 55 & 0.18 & 
				11.47 & 0.51 & 53.61 & 7.27 & 18.44 & 20.16 \\
				& 4 & & 25 & 45490 & 96.7 & 3.3 & 95393 & 166 & 95638 & 55 & 0.26 & 
				26.70 & 0.43 & 51.07 & 9.52 & 18.88 & 20.10 \\
				& 5 & & 25 & 17512 & 91.9 & 8.1 & 99644 & 170 & 99842 & 55 & 0.20 & 
				39.98 & 0.20 & 53.20 & 7.96 & 19.74 & 18.90 \\
				& 6 & & 24.4 & 60674 & 97.0 & 3.0 & 103410 & 148 & 103745 & 54 & 0.33 & 
				67.66 & 0.21 & 52.44 & 10.82 & 17.40 & 19.13 \\
				& 7 & & 25 & 16103 & 89.6 & 10.4 & 98570 & 149 & 98845 & 53 & 0.28 & 
				81.99 & 0.13 & 53.62 & 11.45 & 14.46 & 20.35 \\
				\hline\noalign{\smallskip}
				\multicolumn{3}{c}{\textit{Average}}  & 25 & 25893 & 76 & 24 & 96929 & 159 & 97191 & 56 & 0.27 & 
				33.56 & 0.26 & 52.91 & 7.98 & 19.77 & 19.08 \\
				\hline\noalign{\smallskip}
				\multicolumn{3}{c}{\textit{Total Average}}  & 25 & 26267 & 80 & 20 & 68456 & 126 & 68637 & 44 & 
				0.26 & 17.43 & 0.30 & 53.88 & 9.60 & 15.05 & 21.17 \\
				\hline\noalign{\smallskip}
		\end{longtable}}
	\end{landscape}
}

In Table 3, the next five columns show the contribution of different cost components in the obtained objective function. Columns ``Overtime cost (\%)'', ``Surge capacity cost (\%)'', ``Waiting cost (\%)'', ``Postponement cost (\%)'', and ``OR cost (\%)'' represent the operating rooms’ overtime cost, the total surge cost in downstream units, the patients’ waiting cost, the cost of postponing surgeries, and the fixed cost of opening operating rooms, respectively.

Table 3 shows that the average optimality gap of all instances is 0.26\% that demonstrate the model is working properly and finds near-optimal solutions. Moreover, in the last row of the Table, we observe that the average VSS is 17.43\%, which shows the solution obtained from our proposed stochastic programming model is significantly better than the solution that one would obtain by solving EVP problem. Therefore, we can conclude that considering the uncertainty of random parameters is a critical factor in operating room planning with pooling downstream beds. It is also noteworthy that the VSS increase as the number of specialties increases. In the largest instances with four weeks, VSS raises from 26.70\% to 81.99\% as the number of specialties increases from 4 to 7. The results under columns ``LB Time (\%)'' and ``UB Time (\%)'' also show that the lower bound problem is responsible for the major part of the computational time.

Besides the above analysis, we depicted insightful Figures 5, 6, 7, and 8 to illustrate how the specialties use the shared beds and surging capacities in ICU and wards. These figures help us better understand the dynamics of resource allocation in downstream by our proposed model. We have depicted these figures using the average results on the instances with four weeks of planning horizon and seven specialties. Figure 5 shows the surge capacity in ICU used by different specialties. The sinusoidal trend of the used surge ICU beds is noteworthy; As we can see the cumulative used surge capacity increases during the workdays of the week and then decrease during the weekends. This makes sense because we supposed that our operating room planning problem is related to elective surgeries that are performed during the workdays. It is also very interesting to see that the model forces the specialties to use the surge capacity around the same level even though there is no explicit constraint in this regard in this proposed model.

\afterpage{

\begin{figure}[H]
	\begin{center}
		\hspace{-15pt}\includegraphics[height=2.1in]{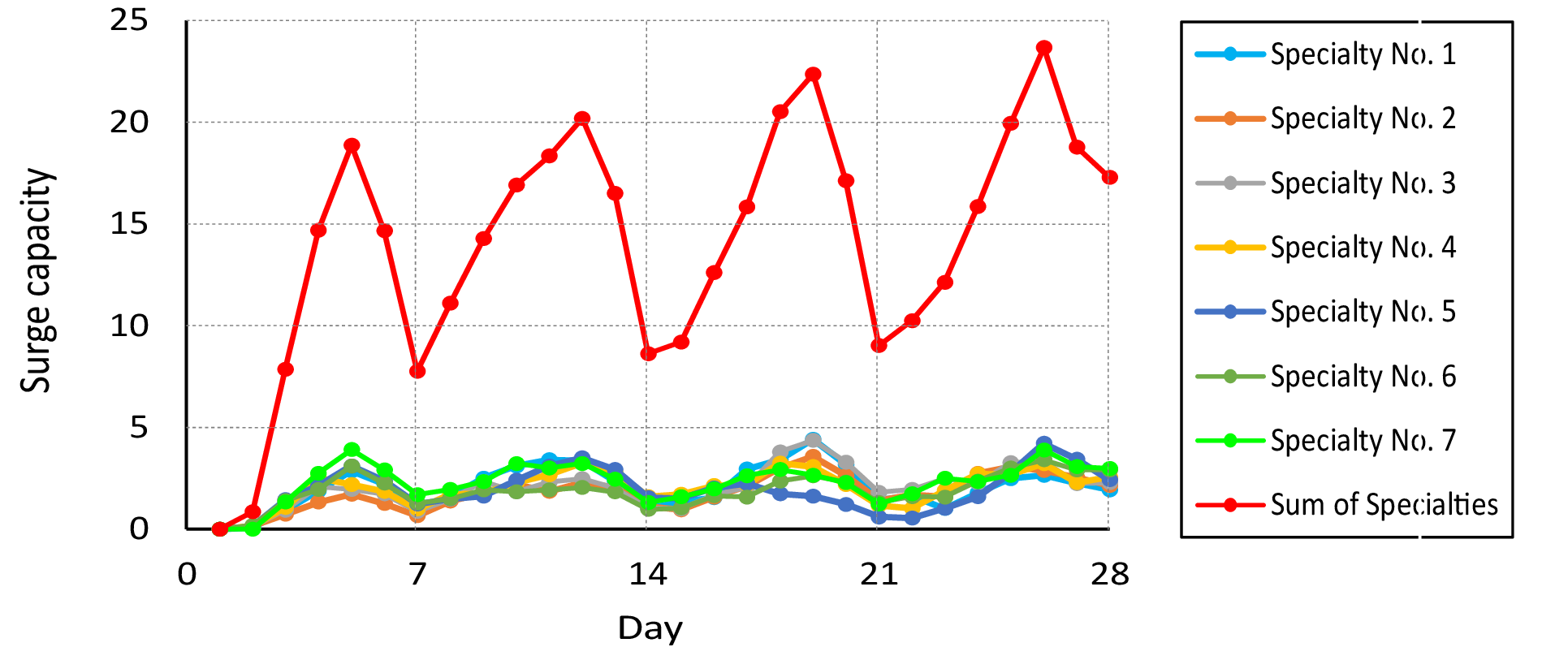}\\[-0pt]
		\caption{Surge capacity used in downstream 1 in the planning horizon.} \label{FigureE.C10-1}
	\end{center}
\end{figure}

\vspace{-40pt}

\begin{figure}[H]
	\begin{center}
		\hspace{-15pt}\includegraphics[height=2.1in]{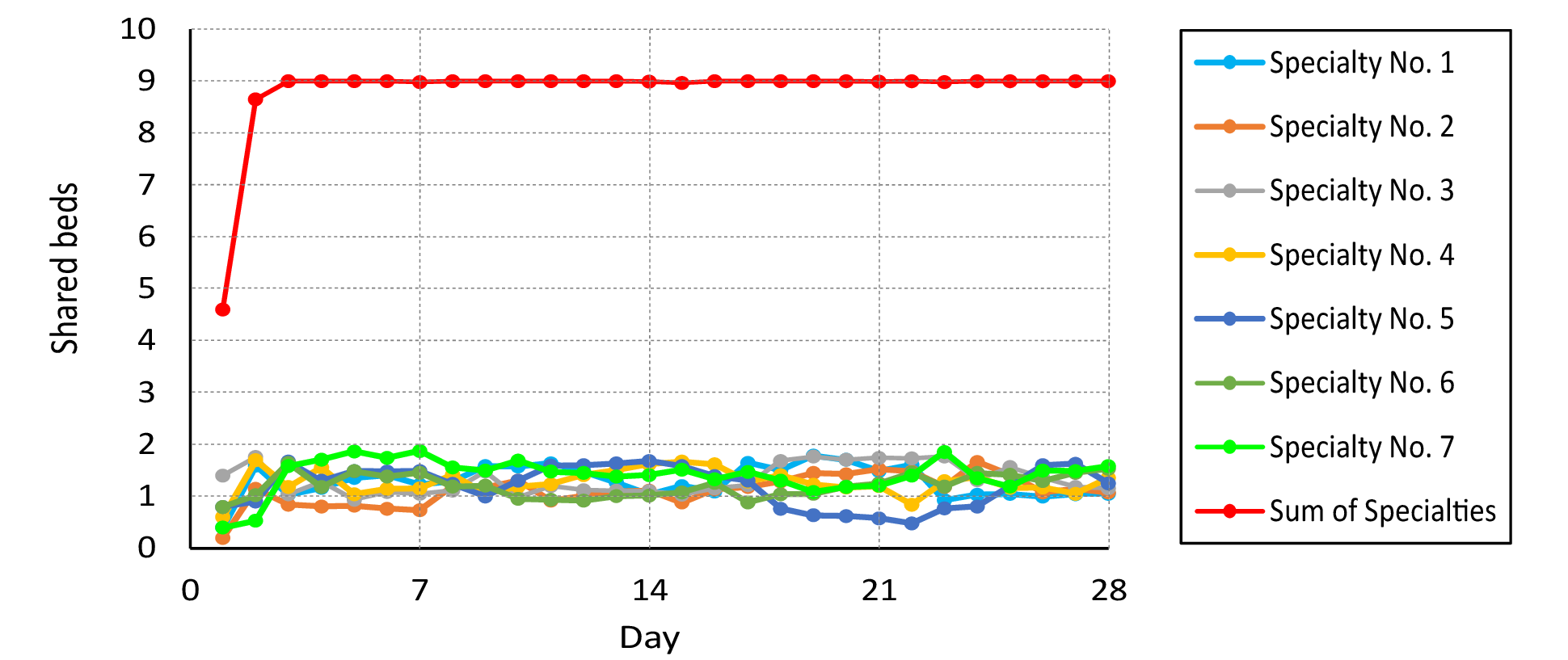}\\[-0pt]
		\caption{The number of shared beds in downstream 1 occupied by patients of different specialties.} \label{FigureE.C10-1}
	\end{center}
\end{figure}

\vspace{-40pt}

\begin{figure}[H]
	\begin{center}
		\hspace{-15pt}\includegraphics[height=2.1in]{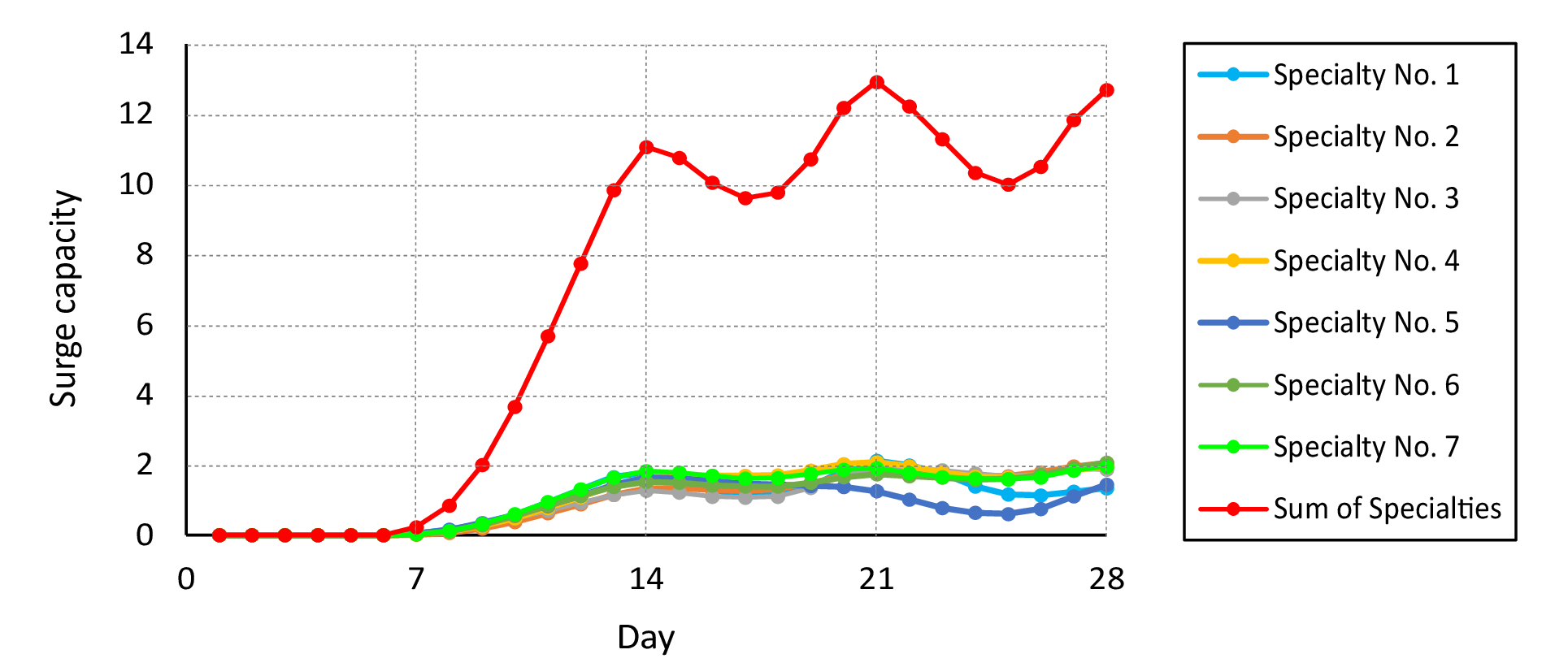}\\[-0pt]
		\caption{Surge capacity used in downstream 2 in the planning horizon.} \label{FigureE.C10-1}
	\end{center}
\end{figure}

\vspace{-40pt}

\begin{figure}[H]
	\begin{center}
		\hspace{-15pt}\includegraphics[height=1.97in]{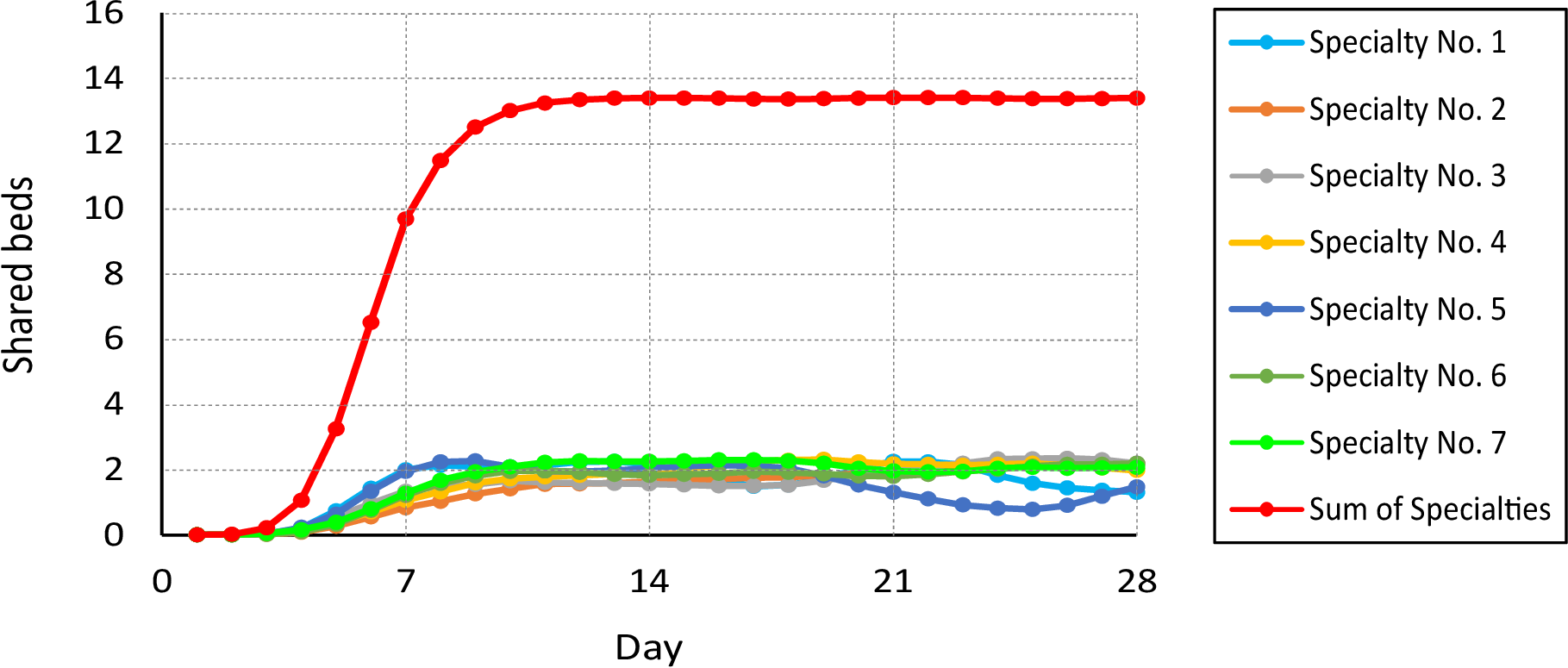}\\[-0pt]
		\caption{The number of shared beds in downstream 2 occupied by patients of different specialties.} \label{FigureE.C10-1}
	\end{center}
\end{figure}

}

Figure 6 depicts the number of shared ICU beds that specialties use during the planning horizon. As we can see in this figure, all shared ICU bed become occupied very quickly at the beginning of the planning horizon after three days. It is also very interesting that all specialties use fairly the same level of shared ICU beds in the planning horizon.

Figure 7 shows the surge capacity used by different specialties in wards during the planning horizon. Compared to Figure 5, the surge capacity in wards increases more slowly and also has less fluctuation. This is because we have more beds in wards and therefore the uncertainty in patients' length of stay in this downstream can be better mitigated using the available resources. 

Figure 8 shows the number of shared beds in wards used by various specialties. We can observe that the cumulative usage increases fairly slowly at the beginning until it reaches the maximum available capacity. Again, we observe that all specialties are fairly using the shared beds during the planning horizon.

\subsection{Sensitivity analysis}

In this section, we perform sensitivity analysis on the costs and uncertain parameters of the proposed model to observe the behavior of the obtained results. In Table 4, we have summarized the results obtained by performing sensitivity analysis on cost parameters independently. In this table, Column ``Parameter'' shows the name of the cost parameter that we have analyzed. Under Column ``Value'', we have different settings for the value of the cost parameters. In our sensitivity analysis, we have modified the values of different cost parameters by $c_{id}^{waiting} := \alpha ^{waiting}c_{id}^{waiting}$, $c^{OR}:=\alpha ^{OR}c^{OR}$, $c_{h}^{bed} := \alpha ^{surge}c_{h}^{bed}$, $c_{i}^{postpone} := \alpha ^{postpone}c_{i}^{postpone}$, and $c^{overtime} := \alpha ^{overtime}c^{overtime}$. The next five columns show the contribution of each cost component to the total cost. ``Waiting Time'' indicator shows the total number of days that all the patients waited before they were operated. ``No. of Postponements'' indicator specifies the number of patients who have been postponed to the next planning horizon. Under ``No. of ORs'', we have the total number of operating rooms that the model has decided to open in the planning horizon.

	\begin{table}[bt]
	%\caption{\hspace{-5pt} Sensitivity analysis of cost parameters.}
	\label{Table1}
	\resizebox{\textwidth}{!}{
		\begin{tabular}{ccccccccccccc}	
			\multicolumn{13}{l}{ \textbf{Table 4} \,\,\, Sensitivity analysis of cost parameters.}\\
			\hline\noalign{\smallskip}		
			 \multicolumn{2}{c}{\textit{Data Info.}} & & \multicolumn{5}{c}{\textit{Cost components}} & & \multicolumn{4}{c}{\textit{Indicators}}\\
			\cline{1-2} \cline{4-8} \cline{10-13} \noalign{\smallskip}		
			 \textit{Parameter} & \textit{Value} & & \textit{\specialcell{Waiting \\[0pt]cost\\[0pt] (\%)}} & \textit{\specialcell{Postpone-\\[0pt]ment \\[0pt]cost\\[0pt] (\%)}} & \textit{\specialcell{OR\\[0pt] cost\\[0pt] (\%)}} & 			 \textit{\specialcell{Overtime \\[0pt] cost \\[0pt] (\%)}} & \textit{\specialcell{Surge\\[0pt] capacity\\[0pt] cost\\[0pt] (\%)}} & & \textit{\specialcell{Waiting\\[0pt] time \\[0pt](Day)}} & 			 \textit{\specialcell{No. of \\[0pt]postpone-\\[0pt]ments}} & \textit{\specialcell{No.\\[0pt] of\\[0pt] ORs}} & \textit{\specialcell{Overtime\\[0pt] (min)}} \\
			\hline\noalign{\smallskip}\\[-14pt]
			\multirow{ 5}{*}{\vspace{17pt} $\alpha ^{waiting}$} & 1 &   & 12.49 & 12.50 & 21.29 & 0.18 & 53.54 & & 79 & 12 & 66.8 & 10 \\[0pt]
			& 3 & & 21.41 & 18.38 & 17.86 & 0.11 & 42.23 & & 44 & 18 & 67.0 & 8 \\[0pt]
			& 5 & & 23.75 & 26.58 & 15.68 & 0.09 & 33.90 & & 31 & 25 & 66.2 & 7 \\[0pt]
			& 7 & & 22.10 & 35.33 & 14.01 & 0.07 & 28.50 & & 23 & 31 & 63.8 & 6 \\[0pt]
			& 10 &   & 16.98 & 47.48 & 12.74 & 0.06 & 22.74 & & 11 & 40 & 62.2 & 5 \\[0pt]
			\hline\noalign{\smallskip}\\[-14pt]
			\multirow{ 5}{*}{\vspace{17pt} $\alpha ^{OR}$} & 1 &   & 12.49 & 12.50 & 21.29 & 0.18 & 53.54 & & 79 & 12 & 66.8 & 10 \\[0pt]
			& 3 & & 11.24 & 12.21 & 40.59 & 0.14 & 35.82 & & 90 & 16 & 59.2 & 11 \\[0pt]
			& 5 & & 10.45 & 12.88 & 49.69 & 0.11 & 26.87 & & 97 & 20 & 54.8 & 11 \\[0pt]
			& 7 & & 10.11 & 14.92 & 54.56 & 0.11 & 20.29 & & 105 & 25 & 51.2 & 13 \\[0pt]
			& 10 &   & 9.30 & 16.63 & 58.89 & 0.11 & 15.07 & & 116 & 32 & 47.4 & 16 \\[0pt]
			\hline\noalign{\smallskip}\\[-14pt]
			\multirow{ 5}{*}{\vspace{17pt} $\alpha ^{surge}$}	& 1 &   & 12.49 & 12.50 & 21.29 & 0.18 & 53.54 & & 79 & 12 & 66.8 & 10 \\[0pt]
			& 3 & & 12.08 & 38.40 & 11.70 & 0.05 & 37.77 & & 105 & 45 & 62.4 & 5 \\[0pt]
			& 5 & & 10.82 & 53.43 & 9.31 & 0.03 & 26.42 & & 97 & 62 & 59.0 & 3 \\[0pt]
			& 7 & & 10.64 & 60.83 & 8.11 & 0.02 & 20.40 & & 93 & 72 & 55.4 & 2 \\[0pt]
			& 10 &   & 9.99 & 67.27 & 7.66 & 0.01 & 15.07 & & 86 & 80 & 54.8 & 1 \\[0pt]
			\hline\noalign{\smallskip}\\[-14pt]
			\multirow{ 5}{*}{\vspace{17pt} $\alpha ^{postpone}$} & 0.5 &   & 9.45 & 29.03 & 22.73 & 0.13 & 38.66 & & 44 & 37 & 62.2 & 6 \\[0pt]
			& 0.75 & & 12.23 & 16.93 & 21.96 & 0.15 & 48.72 & & 71 & 19 & 65.8 & 8 \\[-1pt]
			& 1 & & 12.49 & 12.50 & 21.29 & 0.18 & 53.54 & & 79 & 12 & 66.8 & 10 \\[-1pt]
			& 1.25 & & 13.17 & 5.28 & 21.51 & 0.17 & 59.87 & & 93 & 4 & 68.8 & 10 \\[-1pt]
			& 1.5 &   & 14.24 & 1.99 & 21.95 & 0.15 & 61.66 & & 103 & 1 & 69.8 & 7 \\[-1pt]
			\hline\noalign{\smallskip}\\[-14pt]
			\multirow{ 5}{*}{\vspace{17pt} $\alpha ^{overtime}$} & 1 &   & 12.49 & 12.50 & 21.29 & 0.18 & 53.54 & & 79 & 12 & 66.8 & 10 \\[-1pt]
			& 3 & & 11.82 & 12.70 & 21.51 & 0.43 & 53.54 & & 75 & 12 & 67.6 & 8 \\[-1pt]
			& 5 & & 12.20 & 13.31 & 21.28 & 0.72 & 52.49 & & 78 & 12 & 67.0 & 8 \\[-1pt]
			& 7 & & 12.36 & 13.07 & 21.36 & 0.70 & 52.52 & & 77 & 12 & 67.4 & 6 \\[-1pt]
			& 10 &   & 12.69 & 12.80 & 21.20 & 1.05 & 52.26 & & 82 & 12 & 67.2 & 6 \\[-1pt]
			\hline\noalign{\smallskip}
		\end{tabular}
	}
\end{table}

Table 4 shows that, with an increase in $\alpha^{\text{waiting}}$, the model postpones a higher number of patients. The reason is that, as $\alpha^{\text{waiting}}$ grows, the model seeks to avoid high waiting costs by favoring solutions with lower waiting times. In doing so, the model identifies that postponing some patients and incurring the associated postponement cost is less expensive than treating the same patients while paying a large waiting cost. Therefore, as $\alpha^{\text{waiting}}$ increases, the number of postponed patients also increases. Moreover, the surge capacity cost decreases due to the growth in the number of postponed surgeries. We also observe that increase in $\alpha ^{OR}$ makes the model postpone more patients while the number of opened ORs decreases. As a result, with fewer patients in the system, we have less surge capacity used in downstream. Moreover, Table 4 demonstrates that increasing $\alpha ^{surge}$ results in more postponement of patients to avoid excessive surge capacity cost. 

In Table 4, we set the values of $\alpha ^{postpone}$ to $\{0.5,0.75,1,1.25,1.5\}$ and not $\{1,3,5,7,10\}$ as in other parameter analysis. This is because we realized the values of $\alpha ^{postpone}=\{1,3,5,7,10\}$ does not provide a meaningful analysis since the model does not postpone any surgery for $\alpha ^{postpone}=3$. This trend continued for higher values of $\alpha ^{postpone}$. Therefore, we used $\alpha ^{postpone} \in \{0.5,0.75,1,1.25,1.5\}$. Our results indicate that increasing $\alpha ^{postpone}$ makes the model postpone fewer surgeries as expected. This reduction in the number of postponements results in more scheduled patients and therefore more waiting cost and surge capacity cost. We also observe that, for different values of $\alpha ^{overtime}$, there are no significant changes in cost components and indicators. The only noteworthy point is that for larger values of $\alpha ^{overtime}$ we have less amount of overtime, as expected.

Next, we intend to analyze the sensitivity of the model to the patients' lengths of stay and their surgical durations. To do so, in the revised stochastic programming model, we set $t_{i\omega } := \alpha ^{duration}t_{i\omega } $ and $l_{ih\omega } := \alpha ^{LOS}l_{ih\omega }$. Here, $\alpha ^{duration}$ and $\alpha ^{LOS}$ are the control parameters that we use for sensitivity analysis. Also, as explained in Section 3, $t_{i\omega }$ and $l_{ih\omega }$ denote the patients' surgical durations and lengths of stay, respectively. In our sensitivity analysis, we set $\alpha ^{duration},\alpha ^{LOS} \in \{0.5,0.75,1,1.25,1.5\}$. We have presented the results of both experiments in Table 5.

	\begin{table}[bt]
	%\caption{\hspace{-5pt} Sensitivity analysis of the stochastic parameters.}
	\label{Table1}
	\resizebox{\textwidth}{!}{
		\begin{tabular}{ccccccccccccc}
			\multicolumn{13}{l}{ \textbf{Table 5} \,\,\, Sensitivity analysis of the stochastic parameters.}\\	
			\hline\noalign{\smallskip}		
			\multicolumn{2}{c}{\textit{Data Info.}} & & \multicolumn{5}{c}{\textit{Cost components}} & & \multicolumn{4}{c}{\textit{Indicators}}\\
			\cline{1-2} \cline{4-8} \cline{10-13} \noalign{\smallskip}		
			\textit{Parameter} & \textit{Value} & & \textit{\specialcell{Waiting \\[0pt]cost\\[0pt] (\%)}} & \textit{\specialcell{Postpone-\\[0pt]ment \\[0pt]cost\\[0pt] (\%)}} & \textit{\specialcell{OR\\[0pt] cost\\[0pt] (\%)}} & 			 \textit{\specialcell{Overtime \\[0pt] cost \\[0pt] (\%)}} & \textit{\specialcell{Surge\\[0pt] capacity\\[0pt] cost\\[0pt] (\%)}} & & \textit{\specialcell{Waiting\\[0pt] time \\[0pt](Day)}} & 			 \textit{\specialcell{No. of \\[0pt]postpone-\\[0pt]ments}} & \textit{\specialcell{No.\\[0pt] of\\[0pt] ORs}} & \textit{\specialcell{Overtime\\[0pt] (min)}} \\
			\hline\noalign{\smallskip}\\[-14pt]
			\multirow{ 5}{*}{\vspace{17pt} $\alpha ^{duration}$} & 0.5 &   & 12.35 & 13.40 & 21.41 & 0.00 & 52.85 & & 78 & 12 & 67.0 & 0 \\[0pt]
			& 0.75 & & 12.81 & 11.64 & 21.30 & 0.00 & 54.25 & & 82 & 11 & 66.8 & 0 \\[0pt]
			& 1 & & 12.49 & 12.50 & 21.29 & 0.18 & 53.54 & & 79 & 12 & 66.8 & 10 \\[0pt]
			& 1.25 & & 14.50 & 14.08 & 21.05 & 3.27 & 47.09 & & 96 & 14 & 70.8 & 158 \\[0pt]
			& 1.5 &   & 17.78 & 18.44 & 20.95 & 4.62 & 38.21 & & 121 & 20 & 77.6 & 	307 \\[0pt]
			\hline\noalign{\smallskip}\\[-14pt]
			\multirow{ 5}{*}{\vspace{17pt} $\alpha ^{LOS}$} & 0.5 &   & 31.47 & 0.52 & 54.40 & 0.43 & 13.18 & & 92 & 0 & 70.6 & 10 \\[0pt]
			& 0.75 & & 20.82 & 3.55 & 36.26 & 0.31 & 39.06 & & 85 & 2 & 69.0 & 11 \\[0pt]
			& 1 & & 12.49 & 12.50 & 21.29 & 0.18 & 53.54 & & 79 & 12 & 66.8 & 10 \\[0pt]
			& 1.25 & & 10.94 & 14.93 & 18.56 & 0.13 & 55.44 & & 77 & 16 & 66.4 & 8 \\[0pt]
			& 1.5 &   & 10.23 & 15.55 & 15.35 & 0.11 & 58.76 & & 87 & 19 & 65.4 & 8 \\[0pt]
			\hline\noalign{\smallskip}
		\end{tabular}
	}
\end{table}	

Table 5 indicates that, with increase in $\alpha ^{duration}$, the model gradually increases the number of postponed surgeries and therefore the surge capacity cost decreases. Also, the overtime cost increases significantly due to the increase in surgical times. Moreover, due the limited available times in operating rooms, the model schedules fewer patients than before and therefore the waiting time and its cost increase. 

Moreover, we observe that increase in the patients' lengths of stay results in a significant increase in the surge capacity cost as expected. Also, the model tends to postpone more patients to compromise. This action makes the system less crowded which simultaneously decreases the number of waiting days. An exception in this regard is the case of $\alpha ^{LOS}=1.5$ for which the number of waiting days has increased. This exception is justifiable by observing that the number of opened operating rooms decreases with increases in $\alpha ^{LOS}$ and therefore there is the possibility of having more waiting for some patients.

\section{Conclusion}

In this research, we studied an operating room planning problem with the possibility of pooling the downstream beds among specialties. In our problem, a limited number of beds were available for sharing in ICU and wards. The main sources of uncertainty were the randomness of surgical times and length of stay. We proposed a two-stage stochastic integer programming model and embedded it in a sample average approximation algorithm. In the first stage of our model, the decision maker decides on the allocation of the non-shared beds to specialties, and the allocation of surgeries to operating rooms. Then, in the second stage, after the realization of uncertain parameters, he/she must decide on the allocation of shared ICU and ward beds to specialties during the planning horizon, how many surge capacity beds are required for each specialty in downstream, and also compute the overtime incurred in operating rooms. Our model intended to minimize the total cost including patients' waiting cost and postponement cost, the fixed cost of opening operating rooms, the cost of surge capacity in downstream, and overtime cost. We also proposed a specialized algorithm to enhance the efficiency of the sample average approximation algorithm. 

We performed extensive computational results with three goals: 1) evaluating the effect of three different pooling policies for beds in downstream on the performance of the surgery planning, 2) evaluating the efficiency of the proposed sample average approximation algorithm, and 3) evaluating the behavior and the stability of the proposed model for a large set of instances with various cost and uncertainty parameters.

In the first set of our computational experiments, we compared the results of three downstream-pooling policies, namely, No-Sharing, Midlevel-Sharing, and Full-Sharing policies. In these three strategies, we allow the sharing of 0\%, 50\%, and 100\% of the available beds. Our computational results indicated that the two latter policies lead to 11.29\% and 12.38\% improvement compared to the No-Sharing policy on average, respectively. Also, in some large instances, these improvements were respectively as high as 19.53\% and 16.97\%. In the second part of our computational experiments, we tuned the parameters of our SAA algorithm by depicting the values of optimality gaps, RSD, and solution time for different values of the number of scenarios in lower and upper bound problems and also the number of iterations in SAA.

In the third set of instances, we focused on the Midlevel-Sharing policy and evaluated the efficiency of the proposed SAA algorithm with the tuned parameters. Our results indicated that the proposed approach finds near-optimal solutions with an optimality gap of 0.26\%. Also, the average value of stochastic solution (VSS) is 17.43\%. This value shows that the proposed stochastic programming approach outperforms its corresponding deterministic model by 17.43\%. Moreover, VSS varies between 26.70\% and 81.99\% in instances with four weeks of planning horizon and 4 to 7 specialties. We also provided some insightful figures on the dynamic of shared-beds surge-capacity allocation to specialties. It was very interesting to observe the sinusoidal trend of ICU bed usage and also the fair allocation of shared beds to specialties. We also performed an extensive sensitivity analysis on various parameters and justified the behavior of our model in different settings.

One important limitation of our work is the assumption of full bed sharing among specialties in ICUs and wards. While this policy leads to better system performance, it may be difficult to implement in practice because it requires multi-skilled nursing staff capable of caring for patients from different specialties. This increased flexibility may also raise concerns regarding workload, training requirements, and potential impacts on quality of care.

Future research can focus on an integrated surgery planning and scheduling problem with open scheduling strategies where surgeons are allowed to perform surgeries in multiple operating rooms. Moreover, the concept of pooling can be studied for upstream resources as well as post-anesthesia care units. It is also interesting to study our problem in a robust context and focus on developing a two-stage robust optimization model. Additionally, a promising extension of this work is to incorporate risk-averse decision making by adopting a CVaR-based formulation. Such an approach would allow the planner to explicitly control the risk of extreme outcomes related to downstream congestion, overtime, or excessive surge-capacity usage, which is particularly relevant in high-variability surgical environments.

\vspace{10pt}
\noindent
\textbf{Disclosure of interest}\\
The authors declare that there is no conflict of interest.

\vspace{10pt}
\noindent
\textbf{Funding}\\
This work was supported by the Natural Sciences and Engineering Research Council of Canada (NSERC) under Grant RGPIN-2019-05517.

%\newpage
% Acknowledgments here
%\ACKNOWLEDGMENT{The authors thank the associate editor and two anonymous referees for their constructive comments that significantly improved the quality of this work.}

\vspace{15pt}
\bibliographystyle{ormsv080}
\bibliography{my}

%% Here starts the e-companion (EC)
%%%%%%%%%%%%%%%%%%%%%%%%%%%%%%%%%%%%%%%%%%%%%%%%%%%%%%%%%%
\ECSwitch
\RUNTITLE{Operating room planning with the pooling of downstream beds among specialties}
\numberwithin{equation}{section}
%\ECDisclaimer
%%%%%%%%%%%%%%%%%%%%%%%%%%%%%%%%%%%%%%%%%%%%%%%%%%%%%%%%%%

%%% Main head for the e-companion
\begin{center}
\Large
Electronic Companion 

 ``Operating room planning with the pooling of downstream beds among specialties: A stochastic programming approach''\\
\end{center}

\vspace{10pt}

%\section{Appendix 1- More details on instance generation} \label{Appendix1}
\noindent
\textbf{Appendix 1- More details on instance generation}

\noindent
To generate the surgical time of each patient $i$, we considered normal distribution $N(\mu _{i}^{duration},\sigma _{i}^{duration})$ where $\mu _{i}^{duration}$ and $\sigma _{i}^{duration}$ refer to the average surgical duration of the corresponding specialty and its standard deviation, respectively, and set $\sigma_{i}^{duration}=(1/6)\mu _{i}^{duration}$. The coefficient of (1/6) ensures that the surgical times are generated in $\lbrack \mu_{i}^{duration}-3(\frac{1}{6}\mu _{i}^{duration}),\mu_{i}^{duration}+3(\frac{1}{6}\mu _{i}^{duration})\rbrack $ with a probability of 99.73\%.

Moreover, the last two columns of Table 1 report the average total lengths of stay in hospitals and their standard deviations for patients of different specialties. To randomly generate the length of stay for each patient $i$, we use a normal distribution $N(\mu_{i}^{LOS},\sigma _{i}^{LOS})$. In this distribution, $\mu_{i}^{duration}$ and $\sigma _{i}^{duration}$ respectively refer to the average length of stay for patient $i$ and its standard deviation. To generate a further perturbation in our instances, we supposed that the $\mu _{i}^{LOS}$ of different patients corresponding to the same specialty are not necessarily the same. To address this point, we set $\mu _{i}^{LOS}=[0.75LOS_{s_{i}}^{specialty},1.25LOS_{s_{i}}^{specialty}] $where $LOS_{s_{i}}^{specialty}$ denotes the average length of stay for all patients of the corresponding specialty reported under the third column of Table 1. The patients' lengths of stay generated from $N(\mu_{i}^{LOS},\sigma _{i}^{LOS})$ are corresponding to all downstream units. Therefore, in the end, we divide the generated length of stay $LOS_{i}^{patient}$ between ICU and wards using $LOS_{i}^{ICU}=0.4 LOS_{i}^{patient}$ and $LOS_{i}^{wards}=0.6 LOS_{i}^{patient}$. 

There are also several cost coefficients in our model to be set. We set the fixed cost of opening operating rooms and also the per-minute overtime cost to \$4,437 and \$12.37, respectively (Batun et al. 2011). To set the cost parameters, we consider that each patient $i$ has an urgency level $\alpha_{i}^{priority}$ which are randomly picked from \{1,2,3,4,5\}. Then we set the daily waiting cost of patient $i$ to $\alpha_{i}^{waiting}=\alpha _{i}^{priority}\times 1000.$ Also, the postponement cost of patients is computed by $c_{i}^{postpone}=\alpha _{i}^{priority}\times 15000$. We also set the cost of using one unit of the surge capacity in wards ($c_{1}^{bed}$) to \$62.94. To compute this value, we first divided the annual salary of a wards nurse (\$54,549 as reported on www.ziprecruiter.com) by the total annual work hours (estimated by 52 weeks multiplied by 40 hours) and then divided it by 10, supposing that each nurse takes care of 10 patients on average. The obtained value represents the hourly cost of serving an extra patient in wards. Therefore, we multiplied the recent value by 24 hours and obtained \$62.94 as the daily cost of one unit of surge capacity in wards. Also, for the ICU, we use the same method with the only difference that \$54,549 was replaced by \$95,000 (reported on www.ziprecruiter.com) as the annual salary of an ICU nurse and obtained $c_{2}^{bed}=109.58$ as the daily cost of one extra unit of the surge capacity in ICU.

\vspace{10 pt}

%\bibliographystyle{apalike}
%\bibliography{my}

%%%%%%%%%%%%%%%%%
\end{document}